\documentclass[12pt]{amsart}

\usepackage{amsmath, amssymb, latexsym,amscd}
\begin{document}
\newtheorem{theorem}{Theorem}[section]
\newtheorem{lemma}[theorem]{Lemma}
\newtheorem{proposition}[theorem]{Proposition}
\newtheorem{conjecture}[theorem]{Conjecture}
\newtheorem{condition}[theorem]{Condition}
\newtheorem{example}[theorem]{Example}
\newtheorem{claim}[theorem]{Claim}
\newtheorem{question}[theorem]{Question}
\newtheorem{corollary}[theorem]{Corollary} 
\theoremstyle{definition}
\newtheorem{definition}[theorem]{Definition}
\newtheorem{statement}[theorem]{Statement}
\newtheorem{notation}[theorem]{Notation} 
\theoremstyle{remark}
\newtheorem{remark}[theorem]{Remark}
\newcommand\cl{\begin{claim}}
\newcommand\ecl{\end{claim}}
\newcommand\rem{\begin{remark}\upshape}
\newcommand\erem{\end{remark}}
\newcommand\ex{\begin{example}\upshape}
\newcommand\eex{\end{example}}
\newcommand\nota{\begin{notation}\upshape}
\newcommand\enota{\end{notation}}
\newcommand\dfn{\begin{definition}\upshape}
\newcommand\edfn{\end{definition}}
\newcommand\cor{\begin{corollary}}
\newcommand\ecor{\end{corollary}}
\newcommand\thm{\begin{theorem}}
\newcommand\ethm{\end{theorem}}    
\newcommand\prop{\begin{proposition}}
\newcommand\eprop{\end{proposition}}
\newcommand\lem{\begin{lemma}}
\newcommand\elem{\end{lemma}}
\providecommand\qed{\hfill$\quad\Box$}
\newcommand\pr{{Proof:\;}}
\newcommand\prc{\par\noindent{\em Proof of Claim: }}
\newcommand\dom{{\text{dom}}}
\newcommand\Pn{{P_n}}
\newcommand\ord{{ord}}
\newcommand\rk{{rk}}
\newcommand\car{{\rm{char}}}
\newcommand\M{{{\mathcal M}}}
\newcommand\N{{{\mathcal N}}}
\newcommand\K{{{\mathcal K}}}
\newcommand\D{{{\mathcal D}}}
\newcommand\A{{{\mathcal A}}}
\newcommand\Ps{{{\mathcal P}}}
\newcommand\E{{\mathcal E}}
\newcommand\Y{{{\mathbf{Y}}}}
\renewcommand\O{{{\mathcal O}}}
\renewcommand\L{{{\mathcal L}}}
\newcommand\Th{{\text{Th}}}
\newcommand\IZ{{\mathbb Z}}
\newcommand\IQ{{\mathbb Q}}
\newcommand\IR{{\mathbb R}}
\newcommand\IN{{\mathbb N}}
\newcommand\IC{{\mathbb C}}
\newcommand\F{{\mathcal F}}
\newcommand\IF{{\mathbb F}}
\newcommand\Se{{\mathcal S}}
\newcommand\V{{\mathcal V}}
\newcommand\W{{\mathcal W}}
\newcommand\T{{\mathcal T}}
\newcommand\si{{\sigma}}
\newcommand\n{{\nabla}}
\newcommand\C{{\mathcal C}}
\newcommand\Lr{{{\mathcal L}_{\text{rings}}}}
\newcommand\Lrd{{{\mathcal L}^*_{\text{rings}}}}
\newcommand\B{{\mathcal B}}
\newcommand\La{{\mathcal L}}
\def\U{{ \mathfrak U}}
\def\B{{\mathcal B}}
\def\I{{\mathcal I}}
\def\R{{\mathbb R}}
\def\G{{\mathfrak G}}
\let\le=\leqslant
\let\ge=\geqslant
\let\subset=\subseteq
\let\supset=\supseteq

\newcommand{\bigdoublewedge}{%
  \mathop{
    \mathchoice{\bigwedge\mkern-15mu\bigwedge}
               {\bigwedge\mkern-12.5mu\bigwedge}
               {\bigwedge\mkern-12.5mu\bigwedge}
               {\bigwedge\mkern-11mu\bigwedge}
    }
}

\newcommand{\bigdoublevee}{%
  \mathop{
    \mathchoice{\bigvee\mkern-15mu\bigvee}
               {\bigvee\mkern-12.5mu\bigvee}
               {\bigvee\mkern-12.5mu\bigvee}
               {\bigvee\mkern-11mu\bigvee}
    }
}

 \title{On differential Galois groups of strongly normal extensions}
\date{\today}

 \author{Quentin Brouette}
\email{Quentin.Brouette@umons.ac.be}
  
  \address{D\'epartement de Math\'ematique (Le Pentagone)\\
Universit\'e de Mons\\
20 place du Parc\\ B-7000 Mons\\ Belgium}

  \author{Fran\c coise Point}
  \email{point@math.univ-paris-diderot.fr}
  \address{Research Director at the FRS-FNRS\\D\'epartement de Math\'ematique (Le Pentagone)\\
Universit\'e de Mons\\
20 place du Parc\\ B-7000 Mons\\ Belgium}

\begin{abstract}
 We give a detailed proof of Kolchin's results on differential Galois groups of strongly normal extensions, in the case where the field of constants is not necessarily algebraically closed.  We closely follow former works due to Pillay and his co-authors which were written under the assumption that the field of constant is algebraically closed.
In the present setting, which encompasses the cases of ordered or p-valued differential fields, we find a partial Galois correspondence and we show one cannot expect more in general. 
In the class of ordered differential fields, using elimination of imaginaries in $CODF$, we establish a relative Galois correspondence for definable subgroups of the group of differential order automorphisms.
\end{abstract}
\maketitle

\section{Introduction}
\par  We investigate Galois theory for strongly normal extensions within given classes
of differential fields of characteristic $0$, whose field of constants
are not necessarily algebraically closed.
\par Our main motivation comes from inductive classes $\C$ of topological $\L$-fields
$K$ of characteristic $0$, as introduced in \cite{GP}, and their
expansion to differential structures where the derivation $D$ on $K$
has a priori no interactions with the topology on $K$ (see section 4.1). Let $\C^{ec}$ be the subclass of existentially closed elements of $\C$.
On $\C$, we assume the topology is first-order definable,  the class 
$\C^{ec}$ is elementary and Hypothesis $(I)$ holds. This hypothesis was introduced in \cite{GP} and says that any element of $\C$ has always an extension belonging to $\C$ which contains the field of Laurent series. In particular this implies that the elements of $\C^{ec}$ are large fields. 
Then, one can show that the corresponding class $\C_{D}$ of differential expansions of elements of $\C$ has a model-completion which is axiomatised by expressing that any differential polynomial in one variable with non zero separant, has a zero close to a zero of its associated algebraic polynomial \cite[Definition 5(3)]{GP}.   
\par The setting described above is an adaptation to a topological setting of a former work of M. Tressl \cite{T}. Recently N. Solanki described an alternative approach \cite{Sol} in a similar topological context.

\par Examples of such classes $\C$ are the classes of ordered fields, ordered valued fields, $p$-adic valued fields and valued fields of characteristic zero, and so the corresponding classes $\C^{ec}$ are the classes of real closed fields, real closed valued fields, $p$-adically closed fields and algebraically closed valued fields. When $\C$ is the class of ordered fields, we obtain in this way an axiomatisation of the model completion of $\C_{D}$, the class of closed ordered differential fields, 
alternative to the one given by 
M. Singer \cite{S}; the theory of closed ordered differential fields is denoted by $CODF$. 
\par Given a differential field $F$, we denote by $C_{F}$ its subfield of constants. Let $K\subseteq L$ be two differential fields. Let $\A$ be a saturated differentially closed field containing $L$ and let $\U$ be any intermediate differential field between $L$ and $\A$. 
\par Let $L/K$ be a strongly normal extension as defined by E. Kolchin \cite{K73}, where $C_K$ is not supposed to be algebraically closed. The field $L$ is first assumed to be finitely generated but in the last section we will drop that hypothesis.
 We consider the differential Galois group $gal(L/K)$, namely the group of differential field automorphisms of $L$ fixing $K$, 
and varying $\U$, the differential Galois group $Gal_{\U}(L/K)$, namely the group of differential field automorphisms of $\langle L, C_{\U}\rangle$ fixing $\langle K, C_{\U}\rangle$. When $\U=\A$, we denote this group by $Gal(L/K)$.
\par It follows from the work of Kolchin that $Gal_{\U}(L/K)$ is isomorphic to $H(C_{\U})$, the $C_{\U}$-points of $H$, an algebraic group  defined over $C_{K}$. We give a proof of that fact only using elementary algebra. We also show that there is a bijection between $Gal_{\U}(L/K)$ and $Hom_{K}(L,\U)$, the set of  embeddings of $L$ into $\U$ and this set is in bijection with a definable subset of $\U$ (with parameters in $C_{K}$). 
\par We obtain a partial Galois correspondence. 
Let $E$ be an intermediate extension $(K\subset E\subset L)$ such that $dcl^{\U}(E)\cap L=E$. 
We get that the map $E\mapsto Gal_{\U}(L/E)$ is injective and we observe that it cannot be surjective in general. 
\par We illustrate our results by revisiting the classical examples of strongly normal extensions, namely Picard-Vessiot and Weierstrass extensions.  
\par The existence and unicity of Picard-Vessiot extensions within specific
classes of differential fields (such as formally real or formally
p-adic differential fields) have been considered by several authors
(see for instance \cite{GGO}, \cite{CHP}, using among other things Deligne's work on Tannakian categories). It has been generalised to
strongly normal extensions corresponding to logarithmic differential
equations in \cite{KP}; there, they replace Tannakian categories by the formalism of groupoids as introduced by E. Hrushovski \cite{H}. 
\par In the case of a Picard-Vessiot extension, it is well-known \cite{Kap} that its Galois group is isomorphic to a linear algebraic  group defined over the field of constants. 
\par In the case of a Weierstrass extension of a formally real or a
formally $p$-adic field, we give more explicit information on
which groups may arise as differential Galois groups. When the field of constants is $\IR$ and the differential Galois group $G$ is infinite, we show that it is isomorphic either to $S^{1}$ or $S^{1}\times \IZ/2\IZ$, where $S^{1}$ is the compact connected real Lie group of dimension $1$. When the field of constants is an $\aleph_{1}$-saturated real closed field, we quotient out by $G^{00}$, the smallest type-definable subgroup of $G$ of bounded index \cite{HPP}, and we similarly obtain that $G/G^{00}$ is isomorphic to $S^1$ or to $S^{1}\times \IZ/2\IZ$. When the field of constants is $\IQ_{p}$, we show that $G$ is isomorphic to a subgroup of the $\IQ_{p}$-points of the associated elliptic curve.
\par Assuming now that $T$ is an axiomatisation of a class $\C$ of topological $\L$-fields where Hypothesis $(I)$ holds and letting $T_{c,D}^*$ an axiomatisation of the model-completion of $\C_{D}$, we easily obtain a criterion for the existence of a strongly normal extension within the class of models of $T_{D}$. We compare that result and Theorem 1.5 of \cite{KP}.
\par Then we restrict ourselves to the class of ordered differential fields. 
Let $L/K$ be two differential ordered fields, let $L\subset \U$ with $\U\models CODF$ and let $Aut_{T}(L/K)$ be the subgroup of $Gal_{\U}(L/K)$ of automorphisms that preserve the order. 
Using that the theory $CODF$ has elimination of imaginaries \cite{Point}, we show that if $G$ is a definable subgroup of $H(C_{\U})$, then we can find a tuple $\bar d$ in the real closure of $L$ such that $Aut_{T}(L(\bar d)/K(\bar d))$ is isomorphic to $G\cap H(C_{\U})$. 
\par In the last section we consider non necessarily finitely generated strongly normal extensions $L/K$, as defined by J. Kovacic \cite{Kov1}, and show
that $Gal_{\U}(L/K)$ is isomorphic to a subgroup of an inverse limit of groups of the form $H_i(C_{\U})$, where $H_i$ is an algebraic group defined over $C_{K}$, $i\in I$.
\par This article supersedes the results obtained in \cite{BP} where we dealt with the case of ordered differential fields and are also a part of the first author's thesis \cite[chapter 3]{B}. 
 \section{Preliminaries}\label{basics}
In  this paper, we make use of basic tools coming from model theory in order to show certain definability results for differential Galois groups.
For the convenience of a reader coming from the differential algebra community, we quickly recall some notations coming from model theory and we will be careful pointing out references (mainly from \cite{Ma} and \cite{Sa}).
\medskip
\par Throughout this section, $\L$ denotes a first order language, $\M$ an $\L$-structure and $A\subset M$. Only in this section, we will make the distinction between an $\L$-structure $\M$ and its domain $M$.
\nota We denote by $\L(A):=\L\cup\{c_{a}: a\in A\}$, where $c_{a}$ is a new constant symbol for each element $a\in A$ 
(not to be confused with the element of the subfield of constants in a differential field).  We denote by $Diag(\M)$ the set of all $\L(M)$-sentences true in $\M$. It is a complete $\L(M)$-theory.
\par Given a tuple $\bar a\in M$, we denote by $tp_{A}^{\M}(\bar a)$
the {\it type} of $\bar a$ in $\M$ over $A$.
We will also use the notation $tp_{A}^{T}(\bar a)$, for $T$ a complete theory. 
\par  The set of all algebraic elements over $A$ in $\M$ is denoted by $acl^{\M}(A)$ and called the (model-theoretic) algebraic closure of $A$ in $\M$.
Similarly the set of definable elements over $A$ in $\M$ is denoted by $dcl^{\M}(A)$
 and called the definable closure of $A$ in $\M$.
 \enota
\par We will consider (expansions of) differential fields as $\La_{D}$-structures where $\L$ contains the language of rings $\L_{rings}:=\{+,-,.,0,1\}$ and $\La_{D}:=\La \cup\{^{-1}, D\}$ where $D$ is a new unary function symbol that will be interpreted in our structures as a derivation:
\begin{description}
	\item[Additivity] $\forall a\,\forall b \;\;D(a+b)=D(a)+D(b),$
	\item[Leibnitz rule] $\forall a\,\forall b\;\;\; D(a.b)=a.D(b)+D(a).b .$
\end{description}
\par We denote by $DCF_{0}$ (respectively $ACF_{0}$) the theory of differentially (respectively algebraically) closed fields of characteristic $0$. 
\par A. Robinson showed that $DCF_{0}$ (respectively $ACF_0$) is the model completion of the $\La_{rings,D} \cup\{^{-1}\}$-theory of differential fields (of, respectively, the $\La_{rings} \cup\{^{-1}\}$-theory of fields) of characteristic $0$. Since these theories can be universally axiomatised, it implies that $DCF_{0}$ (respectively $ACF_0$) admits quantifier elimination \cite[Theorem 13.2]{Sa}. 
\par L. Blum gave an elegant axiomatisation of $DCF_0$ \cite[Section 40]{Sa} and showed that any differential field $K$
has an atomic extension, model of $DCF_0$, called the differential closure of $K$ \cite[Section 41]{Sa}. Recall that it is unique up to isomorphism and that the type of a tuple $\bar a$ in the differential closure of $K$ is isolated over $K$; 
following E. Kolchin, 
one says that $\bar a$ is {\it constrained} over $K$ \cite[Section 2]{K74}.  
In the introduction of \cite[page 142]{K74}, 
one can find a "dictionary" between the terminology used by model theorists and by the differential algebra community. 
\medskip
\nota \label{prem} 
Given a field $F$ (respectively a differential field $F$), denote by $\bar F$ the algebraic closure (respectively $\hat{F}$ some differential closure) of $F$. 
We denote by $C_{F}$ the subfield of constants of $F$.
Note that $C_{\hat{F}}=\bar C_{F}$.
\par Given an algebraic group $H$ defined over $F_0\subset F$ and given a subfield $F_0\subset F_1\subset F$, we will denote by $H(F_1)$ its $F_1$-points.
\enota
\par Now let us introduce some differential algebra standard notation.
\nota\label{dfninitial} Let $K\{X_{1},\cdots,X_{n}\}$ be the  ring of differential polynomials over $K$ in $n$ differential indeterminates $X_{1},\cdots, X_{n}$ over $K$, namely it is the ordinary polynomial ring in indeterminates $X_{i}^{(j)}$, $1\leq i\leq n$, $j\in \omega$, with by convention $X_{i}^{(0)}:=X_{i}$. One can extend the derivation $D$ of $K$ to this ring by setting $D(X_{i}^{(j)}):=X_{i}^{(j+1)}$ and using the additivity and the Leibnitz rule. 

\par Set $\bar X:=(X_{1},\cdots, X_{n})$ and $\bar X^{(j)}:=(X_{1}^{(j)},\cdots, X_{n}^{(j)})$. 
We will denote by $K\langle \bar X\rangle$ the fraction field of $K\{\bar X\}$.
Let $f(\bar X)\in K\{\bar X\}\setminus K$ and suppose that $f$ is of order $m$, then we associate with 
$f(\bar X)$ the ordinary polynomial $f^*(\bar X^*)\in K[\bar X^*]$ with $\bar X^*$ a tuple of variables of the appropriate length such that $f(\bar X)=f^*(\bar X^{\n})$ where $\bar X^{\n}:=(\bar X^{(0)},\ldots ,\bar X^{(m)})$.
 We will make the following abuse of notation: if $\bar b\in K^n$, then $f^*(\bar b)$ means that we evaluate the polynomial $f^*$ at the tuple $\bar b^{\n}$.

\par If $n=1$, recall that the separant $s_{f}$ of $f$ is defined as $s_{f}:=\frac{\partial f }{\partial X_{1}^{(m)}}$.
\enota
\medskip
\par We will use the following consequence of the quantifier elimination results for $T=ACF_{0}$ (respectively $DCF_{0}$). 
Let $\M\models T$, let $K$ be a subfield of $M$ (respectively a differential subfield) and $\bar a\in M$, 
then the type $tp_{K}^{T}(\bar a)$ of $\bar a$ over $K$ is determined by the ideal $$\I_{K}(\bar a):=\{p(\bar X)\in K[\bar X]:\;p(\bar a)=0\}$$ (respectively the differential ideal $\I^D_{K}(\bar a):=\{p(\bar X)\in K\{\bar X\}:\;p(\bar a)=0\}$).

\par Finally letting $F_{0}\subset F$, with $F\models ACF_{0}$ and $F_{0}$ a subfield, recall that given an $F_{0}$-definable function $f$ in $F$, we can find finitely many $F_{0}$-algebraic subsets $X_{i}$ such that $f$ is equal to an $F_{0}$-rational function on $X_{i}$ (straightforward from \cite[Proposition 3.2.14]{Ma}).  
Moreover if $E$ is an $F_0$-definable equivalence relation on $F^{n}$, then there is an $F_{0}$-definable function $f$ from $F^{n}$ to $F^{\ell}$ 
such that $E(\bar x,\bar y)$ iff $f(\bar x)=f(\bar y)$ \cite[Theorem 3.2.20]{Ma}. This is another way to express that $ACF_0$ admits elimination of imaginaries (e.i.)
(see \cite{M}, Chapter I, section 4, for a formal definition of e.i.)
\medskip

\par We will work inside a {\it saturated} structure $\M$ \cite[Definition 4.3.1]{Ma}, but as usual we will only need $\kappa^+$-saturation for some cardinal $\kappa$.
Such a structure $\M$ has the following {\it homogeneity} property \cite[Proposition 4.3.3 and Proposition 4.2.13]{Ma}.
Let $A\subset M$ such that $|A|<|M|$ and $\bar a, \bar b$ be finite tuples of elements of $M$ 
with $tp^{\M}_A(\bar a)=tp^{\M}_A(\bar b)$, 
then there is an automorphism of $\M$ fixing $A$ and sending $\bar a$ to $\bar b$.
\medskip

\section{Strongly normal extensions}

\par Let $K\subseteq L\subseteq \U$ be three differential fields and let $C_K\subset C_L\subset C_{\U}$ their respective subfields of constants; we will also use the notations $\U/L/K$, $L/K$ or $\U/L$.  Let $\A$ be a saturated model of $DCF_0$ containing $\U$; $\A$ will play the role of a universal extension \cite[page 768]{K53}. 
\par Since $DCF_{0}$ admits quantifier elimination, $\A$ has the following property: any $\L_{rings,D}$-embedding from $K\langle \bar a \rangle$ to $K\langle \bar b \rangle $ fixing $K$ and sending $\bar a$ to $\bar b$, where $\bar a, \bar b$ are finite tuples from $\A$, can be extended to an automorphism of $\A$ fixing $K$ (we use the homogeneity of $\A$, see section \ref{basics}). 
\par Denote by $Hom_{K}(L,\A)$ the set of $\L_{rings,D}$-embeddings from $L$ to $\A$ which are fixed on $K$. 
\medskip

\par The notion of {\it strongly normal} extension has been introduced by E. Kolchin.  
Recall that an element $\tau\in Hom_{K}(L,\A)$ is {\it strong} if $\tau$ is the identity on $C_{L}$ and if $\langle L,C_{\A} \rangle=\langle\tau(L),C_{\A} \rangle$  \cite[page 389]{K73}.
Then $L$ is a {\it strongly normal extension} of $K$ if $L$ a finitely generated extension of $K$ such that every $\tau\in Hom_{K}(L,\A)$ is {\it strong} \cite[page 393]{K73}. 
\par From the above remarks, it follows that, whenever $L$ is a finitely generated extension of $K$, any element of $Hom_{K}(L,\A)$ extends to an $\L_{rings,D}$-automorphism of $\A$ fixing $K$ and when $L$ is a strongly normal extension of $K$, any such automorphism restricts to an automorphism of $\langle L,C_{\A} \rangle.$
\par One can extend this notion to not necessarily finitely generated extensions $L$ of $K$ by asking that $L$ is the union of finitely generated strongly normal extensions \cite[Chapter 2, section 1]{Kov1}. Except in Section \ref{nfg}, we will assume that strongly normal extensions are finitely generated.
\par By \cite[Chapter 6, section 3, Proposition 9]{K73}, if $L$ is a strongly normal extension of $K$, then $C_{K}=C_{L}$. 
\nota \label{auto-nota} Let $\U/L/K$ be two differential extensions of $K$. We denote by $gal(L/K)$ the group of $\L_{rings,D}$-automorphisms of $L$ fixing $K$, by $Gal(L/K)$ the group of $\L_{rings,D}$-automorphisms of $\langle L, C_{\A}\rangle$ fixing $\langle K, C_{\A}\rangle$, and by $Gal_{\U}(L/K)$ the group of $\L_{rings,D}$-automorphisms of $\langle L, C_{\U}\rangle$ fixing $\langle K, C_{\U}\rangle$.
\enota
In the following definition, we use a similar terminology as used by I. Kaplansky for a differential extension to be normal \cite[page 20]{Kap}.
\dfn\label{normal} Let $L$ be a differential field extension of $K$ with $K\subset L\subset \A$. Let $\Sigma\subset Hom_{K}(L,\A)$. Then $L$ is $\Sigma$-normal if for any $x\in L\setminus K$, there is $\tau\in \Sigma$ such that $\tau(x)\neq x$.
\edfn
\par I. Kaplansky showed that if $K$ is a differential field of characteristic $0$ with algebraically closed field of constants, then any Picard-Vessiot extension $L$ of $K$ is $gal(L/K)$-normal \cite[Theorem 5.7]{Kap}. 
E. Kolchin showed that if $L$ is a strongly normal extension of $K$, then $L$ is $Hom_{K}(L,\A)$-normal \cite[Chapter 6, section 4, Theorem 3]{K73}.
\medskip
\par Let $L$ be a strongly normal extension of $K$, then E. Kolchin has shown that $Gal(L/K)$ is isomorphic to the $C_{\A}$-points of an algebraic group $H$ defined over $C_{K}$, which is equivalent, since $ACF_{0}$ admits elimination of imaginaries  \cite[Chapter 16]{Poizat}, to prove that it is isomorphic to a group interpretable in $C_{\A}$ (with parameters in $C_K$). 
\par In the following, we will revisit Kolchin's result and give a detailed proof that there is an algebraic group $H$ defined over $C_{K}$ such that $Gal_{\U}(L/K)$ is isomorphic to the $C_{\U}$-points of $H$ for any field $L\subset \U\subset \A$ ( Theorem \ref{Kolchin}). 
We also show that $H(C_{\U})$ is in bijection with the following sets:
$Hom_{\langle K, C_{\U}\rangle}(L,\U)$, $Hom_{K}(L,\U)$ and a definable subset of $\U$ with parameters in $C_{K}$.
\par We begin by a few easy remarks. 
\rem\label{ld} Suppose that $C_{K}$ is algebraically closed in $C_{\U}$ (denoted by $C_K\subset_{ac} C_{\U}$). 
Then the fields $\U$ and $\bar C_{K}$ are linearly disjoint over $C_{K}$, respectively $\U$ and $\langle K,\bar C_{K}\rangle$ are linearly disjoint over $K$. Therefore $\langle L,\bar C_{K}\rangle \cap \U=L$.
\erem
\pr  
By assumption $C_{K}$ is algebraically closed in $C_{\U}$ and so in $\U$. Note that since $\U$ is separable over $C_{K}$, $\U$ is a regular extension of $C_{K}$. Then the first assertion follows by \cite[Chapter 3, Theorem 2]{Lang} and the second one from \cite[Chapter 3, Section 1, Proposition 1]{Lang}.  
\par Let $u:=\frac{\sum_{i}\ell_{i}.c_{i}}{\sum_{j}\ell_{j}.c_{j}}\in \langle L,\bar C_{K}\rangle \cap \U$ with $\ell_{i},\ell_{j}\in L$, $c_{i}, c_{j}\in \bar C_{K}$. Since $\bar C_{K}$ and $\U$ are linearly disjoint over $C_{K}$, there exist $d_{i}, d_{j}\in  C_{K}$ such that $u.\sum_{j}\ell_{j}.d_{j}=\sum_{i}\ell_{i}.d_{i}$, so $u\in L$.
\qed
\rem\label{ldA} Note that we have always that $\U$ (respectively $L$) and $C_{\A}$ are linearly disjoint over $C_{\U}$ (respectively $C_{K}$) and so $\langle L, C_{\A}\rangle \cap \U=\langle L,C_{\U}\rangle$.
\erem
\pr Indeed, $C_{\U}\subset_{ac} \U$ and $DCF_{0}$ is model-complete.\qed
\rem\label{star} Assume now that $L/K$ is a strongly normal extension. 
Then for any $\si\in Hom_{K}(L,\U)$, we have $\langle L, C_{\U}\rangle =\langle \si(L), C_{\U}\rangle$. 
\erem
\pr  
Given $\si\in Hom_{K}(L,\U)$, since $L$ is a strongly normal  extension of $K$, we have $\langle L, C_{\A}\rangle =\langle \si(L), C_{\A}\rangle$.
\par By Remark \ref{ldA}, we have that $\langle L, C_{\A}\rangle\cap \U=\langle L, C_{\U}\rangle$. Since $\si(L)\subset \U$, applying the same remark to $\si(L)$, we also have that $\langle \si(L), C_{\A}\rangle\cap \U=\langle \si(L), C_{\U}\rangle$. 
Therefore, we obtain: $\langle L, C_{\U}\rangle =\langle \si(L), C_{\U}\rangle$. 
\qed
\medskip
\rem\label{L} Let $\bar u\in L$, then $tp_{K}^{\A}(\bar u)\models tp_{\langle K,\bar C_{K}\rangle}^{\A}(\bar u)$.
\erem
\pr 
Since $C_L=C_K$, $L$ and $\bar C_{K}$ are linearly disjoint over $C_{K}$, which implies that $L$ and $\langle K,\bar C_{K}\rangle$ are linearly disjoint over $K$. So, the ideal $\I^D_{\langle K,\bar C_{K}\rangle}(\bar u)$ (Notation \ref{prem}) has a set of generators in $K\{\bar X\}$ \cite[Chapter 3, Section 2, Theorem 8]{Lang}.
Since the theory $DCF_{0}$ admits q.e.,  the type of a tuple is determined by the set of differential polynomials it annihilates. 
\qed
\medskip
\par It is well-known that a strongly normal extension $L$ of $K$ is of the form $K\langle \bar a\rangle$ for some tuple $\bar a\in \hat{K}$ \cite[section 9, Theorem 3]{K74}, which implies that the definition of {\it strongly normal} does not depend on the universal extension one is working in. 
\par Moreover one can observe that ''$K\langle \bar a\rangle$ is a strongly normal extension of $K$'' only depends on $tp_{K}^{DCF_{0}}(\bar a)$. This is well-known but we will use the proof later and so we will quickly review it below, following \cite{P98}. We derive the additional information that, in case $C_{K}\subseteq_{ac} C_{\U}$, any tuple $\bar b\in \hat K\cap \U$ in the image of some $\si\in Hom_{K}(L,\A)$, belongs to $L$.
\prop\label{Tsn} {\rm (See \cite[section 9, Theorem 3]{K74}, \cite{P98})} Let $L/K$ be a strongly extension. Then 
$L$ is of the form $K\langle \bar a \rangle $ for some constrained tuple $\bar a$ over $K$. Moreover there exists a formula $\xi(\bar y)\in tp_{K}^{DCF_{0}}(\bar a)$ such that for any $\bar b\in \A$ with $\A\models \xi(\bar b)$, we have $\bar a\in \langle K, \bar b,C_{\A} \rangle$, $\bar b\in \langle K, \bar a,C_{\A} \rangle$ and $C_{\langle K, \bar b\rangle}=C_{K}$. 
\par Moreover, assume that $L\subseteq \U$ and $\bar b\in \U$. Then if $\U\models \xi(\bar b)$, we have $\bar b\in\langle L,C_{\U}\rangle)$ and if in addition, $C_{K}\subseteq_{ac} C_{\U}$ and $\bar b\in\U\cap \hat{K}$, then  
$\bar b\in L$.
\eprop
\pr Let $\bar a:=(a_{1},\cdots,a_{n})\in L$ be such that $L:=K\langle \bar a \rangle$ and let $\Sigma(\bar y)= tp_{K}^{DCF_{0}}(\bar a)$. Let 
$$\chi_{i\eta}(a_{i},\bar y):=\exists \bar c_{1i\eta}\exists \bar c_{2i\eta}\;(a_{i}=\frac{p_{1i\eta}(\bar y,\bar c_{1i\eta})}{p_{2i\eta}(\bar y,\bar c_{2i\eta})}\;\&\;p_{2i\eta}(\bar y,\bar c_{2i\eta})\neq 0\;\&\;D(\bar c_{1i\eta})=D(\bar c_{2i\eta})=0),$$
where $p_{1i\eta}(\bar Y,\bar Z),\;p_{2i\eta}(\bar Y,\bar Z)$, $1\leq i\leq n,$ range over all elements of $K\{\bar Y\}[\bar Z]\setminus\{0\}$, $\bar Z$ is a countably infinite tuple of variables. $(\star)$
\cl \label{comp} $\Xi(\bar y):=\Sigma(\bar y)\cup Diag(L)\cup\{\bigvee_{i=1}^n\neg\chi_{i\eta}(a_{i},\bar y); \eta<\vert K\vert\}$ is inconsistent.
\ecl
 \prc Otherwise, since $\A$ is $\vert K\vert^+$-saturated, $\Xi(\bar y)$ would have a realisation 
$\bar b\in\A$, realising the same type in $\A$ over $K$ as $\bar a$. So there would exist an elementary map from $K\langle \bar a \rangle $ to $K\langle \bar b \rangle $ fixing $K$ and sending $\bar a$ to $\bar b$.  
However, $\bar a$ would not belong to $\langle K\langle \bar b\rangle,C_{\A}\rangle$ contradicting the fact that $L$ is a strongly normal  extension of $K$. 
\qed
\medskip
\par So by Claim \ref{comp} and the compactness theorem \cite[Theorem 2.1.4]{Ma}, there exists a formula $\phi(\bar y)$ in $tp_{K}^{DCF_{0}}(\bar a)$ such that $DCF_{0}\cup Diag(L)\models \forall \bar y (\phi(\bar y)\rightarrow (\bigvee_{j\in J}\bigwedge_{i=1}^n \chi_{ij}(a_{i},\bar y))$, with $J$ finite.
So, $\hat{L}\models \forall \bar y (\phi(\bar y)\rightarrow (\bigvee_{j\in J}\bigwedge_{i=1}^n \chi_{ij}(a_{i},\bar y)))$.
\par We may assume that $\hat{K}\preceq \hat{L}\preceq \A$ and so there exists $\bar b_{1}\in \hat{K}$ such that $\hat{K}\models  \phi(\bar b_{1})$ and so $\hat{L}\models \phi(\bar b_{1})$. Since $C_{\hat{L}}=\bar C_{L}=\bar C_{K}$, we get $\langle K, \bar b_{1},C_{\hat{L}}\rangle\subset \hat{K}$, and so $\bar a\in \hat{K}$. 
\par So w.l.o.g. we may choose a quantifier-free formula $\xi(\bar y)$ isolating the type of $\bar a$ over $K$.
 Since $L$ is a strongly normal extension of $K$, we also have that 
$DCF_{0}\cup Diag(L)\models  \forall \bar y\;(\xi(\bar y)\rightarrow 
(\bigvee_{j'\in J'}\bigwedge_{i=1}^n\,\chi_{ij'}(y_{i},\bar a)))$, with $J'$ finite. Therefore, if $\bar b\in \A$ and $\xi(\bar b)$ holds, then we have $\bar a\in \langle K, \bar b,C_{\A} \rangle$, $\bar b\in \langle K, \bar a,C_{\A} \rangle$.
\par Finally note that the formula $\xi(\bar y)$ implies that $D(f(\bar y))\neq 0$, for any $f(\bar Y)\in K\langle \bar Y\rangle$ such that $f(\bar a)$ is well-defined and $f(\bar a)\notin K$. So, if $\A\models \xi(\bar b)$, then $K\langle \bar b\rangle/K$ is a strongly normal  extension.
\par If $\bar b\in \U$, then $\bar b\in \langle L, C_{\A} \rangle\cap \U$, which is equal to $\langle L, C_{\U} \rangle$ (Remark \ref{ldA}).
\par If in addition $C_K\subset_{ac} C_{\U}$, for any $\bar b\in \U\cap \hat K$, since $\bar C_{L}=\bar C_{K}=C_{\hat{K}}$, we have that $\bar b\in \langle L, \bar C_{K} \rangle\cap \U$. By Remark \ref{ld}, $\bar b\in L$.
This finishes the proof of the proposition. \qed
\bigskip
\par Note that by a recent result of Pogudin \cite{Pog}, any finitely generated differentially algebraic field extension $L$ of $K$ such that $C_L\neq L$, is generated by one element; this extends in characteristic $0$ a former result of Kolchin \cite[Chapter 2, Proposition 9]{K73}. Therefore we may assume that any strongly normal extension $L/K$ is of the form $K\langle a\rangle$ for a single element $a$.
\rem\label{uni} 
\par Keeping the same notations as in  the proof above but assuming that $n=1$ ($a_i=a$) and using the terms occurring in $(\star)$, we will check that we can express $a$ as a rational function of $b$, a fixed tuple $\bar k\in K$ and tuples $\bar c$ either belonging to $C_{\A}$ or to $C_{\U}$ or $C_{K}$ and conversely with $b$ in place of $a$. 
We will use the fact that we may use the same rational functions later, whenever $\bar c\in C_{\A}, C_{\U}$ or $C_{K}$ (see Theorems \ref{Gal}, \ref{Kolchin}).
\par Indeed, in the above proposition, we have seen that for $b\in \U$, if $\xi(b)$ holds, then $b\in \langle L,C_{\U}\rangle$ and further for $b\in \U\cap \hat{K}$ and if $C_K\subset_{ac} C_{\U}$, then $b\in L$. 
 \par $(1)$ First take $b\in \A$ with $\A\models \xi(b)$,
rewrite the polynomials $p(Y,\bar c)$ occurring in $\chi_{j},\;\chi_{j'}$, with  $p(Y,\bar c)\in K[\bar c]\{Y\}$,\;$\bar c\in C_{\A}$, as $\tilde p(Y,\bar T,\bar c)=\sum_{s} c_{s}m_{s}(Y,\bar T)$ with $c_{s}\in C_{\A}$ and $m_{s}(Y,\bar T)$ a monomial in $Y^{\n}$ and $\bar T$.
We obtain $\tilde p_{2}(b,\bar k,\bar c_{2}).a-\tilde p_{1}(b,\bar k,\bar c_{1})=0$ and $\tilde p_{2}(b,\bar k,\bar c_{2})\neq 0$, with $\tilde p_{2}(Y,\bar T,\bar Z_{2}).X=\sum_{s} Z_{2s}m_{s}(Y,\bar T).X$ and $\tilde p_{1}(Y,\bar T,\bar Z_{1})=\sum_{s} Z_{1s}m_{s}(Y,\bar T)$. Set $\bar \alpha_{1}:=(c_{1s})$, $\bar \alpha_{2}:=(c_{2s})$ and $\bar \alpha:=(\bar \alpha_{1},\bar \alpha_{2})$ and $\bar Z:=(\bar Z_{1},\bar Z_{2})$ with $\bar Z_{1}$, respectively $\bar Z_{2}$, of the same length as $\bar \alpha_{1}$, respectively $\bar \alpha_{2}$. 
Denote by $t_{1}(Y,\bar T,\bar Z_{1}), t_{2}(Y,\bar T,\bar Z_{2})\in \IZ\{Y\}[\bar T,\bar Z]$ the terms corresponding to $\sum_{s} Z_{1s}.m_{s}(Y,\bar T)$, $\sum_{s} Z_{2s}m_{s}(Y,\bar T).X$. We have that $$\sum_{s} c_{1s}.m_{s}(b,\bar k)-\sum_{s} c_{2s}m_{s}(b,\bar k).a=0\;\;(\star).$$
Let $f(Y,\bar T,\bar Z):=\frac{t_{1}(Y,\bar T,\bar Z_{1})}{t_{2}(Y,\bar T,\bar Z_{2})}$. Similarly we can express $b$ using rational functions $\tilde f:=\frac{s_{1}(Y,\bar T,\bar Z_{1})}{s_{2}(Y,\bar T,\bar Z_{2})}\in \IQ\langle Y,\bar T,\bar Z\rangle$ of $a,\bar k,\bar \alpha$, namely as $b=\tilde f(a, \bar k, \bar \alpha)$. These rational functions vary in a finite set $\{f_{\ell},\;\tilde f_{\ell}; \ell\in \Lambda\}$. W.l.o.g. we enlarge the tuples $\bar k\in K$, $\bar \alpha$ with $\bar \alpha\in C_{\A}$, in order to use them in both expressions (when we express $a$ in terms of $b$ and $b$ in terms of $a$).
\par $(2)$ Assume in addition that $b\in \U$. Choose a basis $(e_{\mu})$ of $C_{\A}$ over $C_{\U}$ and express each $c_{1s}\in C_{\A}$ as $\sum_{\mu} \alpha_{1s\mu}.e_{\mu}$ with $\alpha_{1s\mu}\in C_{\U}$ and similarly for  $c_{2s}$. If we substitute it in expression $(\star)$, we obtain for all $\mu$ that:
$\sum_{s} \alpha_{1s\mu}.m_{s}(b,\bar k)-\sum_{s} \alpha_{2s\mu}.m_{s}(b,\bar k).a=0$, since $\U$ and $C_{\A}$ are linearly disjoint over $C_{\U}$. \par Moreover, for some $\mu_{0}$ we have that $\sum_{s} \alpha_{2s\mu_{0}}.m_{s}(b,\bar k)\neq 0$.
Set $\bar \alpha_{1}:=(\alpha_{1s})_{s}$, $\bar \alpha_{2}:=(\alpha_{2s\mu})_{s}$ and $\bar \alpha:=(\bar \alpha_{1},\bar \alpha_{2})$. 
\par We obtain for $b\in \U$ that $a=f_{\ell}(b,\bar k,\bar \alpha)$ and $b=\tilde f_{\ell}(a,\bar k,\bar \alpha)$, with $\bar \alpha$ varying in $C_{\U}$. This implies in particular that $b\in \langle L,C_{\U}\rangle$.
\par $(3)$ If we restrict ourselves to $b\in L$, we get that $a=f_{\ell}(b,\bar k,\bar \alpha)$ and $b=\tilde f_{\ell}(a,\bar k,\bar \alpha)$, with $\bar \alpha$ varying in $C_{K}$. (We proceed as above choosing now a basis $(e_{\mu})$ of $\bar C_{K}$ over $C_{K}$ and using that $L$ and $\bar C_{K}$ are linearly disjoint.) 
\par $(4)$ Finally assume that $b\in \hat{K}\cap \U$ and in addition that $C_{K}\subset_{ac} C_{\U}$. 
Since now $\U$ and $\bar C_{K}$ are linearly disjoint over $C_{K}$, 
we may assume that we have chosen a basis $(e_{\mu})$ of $C_{\A}$ over $C_{\U}$ which extends a basis of $\bar C_{K}$ over $C_{K}$. We express each $c_{1s}\in \bar C_{K}$ as $\sum_{\mu} \alpha_{1s\mu}.e_{\mu}$ with $\alpha_{1s\mu}\in C_{K}$ (and similarly for $c_{2s}$).
\par We obtain that $a=f_{\ell}(b,\bar k,\bar \alpha)$, respectively $b=\tilde f_{\ell}(a,\bar k,\bar \alpha)$, with $\bar \alpha\in C_{K}$, which implies that $b\in L$.
\erem

\lem\label{type} Let $\bar a\in \A$ be such that $tp^{\A}_{\langle K,C_{\bar K}\rangle}(\bar a)$ is isolated, then $tp_{\langle K,C_{\bar K}\rangle}^{\A}(\bar a)\models tp_{\langle K, C_{\A}\rangle}^{\A}(\bar a)$.
 \elem
\pr Let $\phi(\bar y)$ be an $\L_{rings,D}(\langle K,C_{\bar K}\rangle)$-formula isolating $tp_{\langle K,C_{\bar K}\rangle}^{\A}(\bar a)$. 
By the way of contradiction, suppose there is a tuple $\bar a_{1}\in \A$ satisfying the set of formulas $tp_{\langle K,C_{\bar K}\rangle}^{\A}(\bar a)\cup\{\neg \psi(\bar y)\}$ where $\psi(\bar y)\in tp_{\langle K, C_{\A}\rangle}^{\A}(\bar a)$. Rewrite $\psi(\bar y)$ as $\tilde \psi(\bar y,\bar c)$, a $\L_{rings, D}(K)$-formula:  where $\bar c\in C_{\A}$. So, we have: 
$$\A\models \exists \bar x_{1}\exists \bar x_{2}\exists \bar c\;(\phi(\bar x_{1})\;\&\;\phi(\bar x_{2})\;\&\;D(\bar c)=0\;\&\;\tilde \psi(\bar x_{1},\bar c)\,\&\,\neg \tilde \psi(\bar x_{2},\bar c)).$$ Therefore $\widehat{K}\models \exists \bar x_{1}\exists \bar x_{2}\exists \bar c\;(\phi(\bar x_{1})\;\&\;\phi(\bar x_{2})\;\&\;D(\bar c)=0\;\&\;\tilde\psi(\bar x_{1},\bar c)\,\&\,\neg \tilde\psi(\bar x_{2},\bar c))$,
a contradiction.
\qed

\lem\label{emb}  Let $L/K$ be a strongly normal extension, then there is an embedding of $gal(L/K)$ into $Gal(L/K)$. 
Any element of $Hom_{K} (L,\U)$ extends uniquely to an element of $Gal(L/K)$, which restricts to an element of $Gal_{\U}(L/K)$.
\elem
\pr Let $a\in L$ be such that $L=K\langle a \rangle$.
\par Let $\si\in Hom_{K} (L,\U)$. Since $DCF_{0}$ admits q.e., $tp_{K}^{\A}(a)=tp_{K}^{\A}(\si(a))$ and these types are determined by $\I^D_{K}(a)$ (respectively $\I^D_{K}(\si(a))$).
Since $a\in L$, we have that $tp_{K}^{\A}(a)\models tp_{\langle K, \bar C_{K}\rangle}^{\A}(a)$ (Remark \ref{L}).
Since $tp_{K}^{\A}(a)=tp_{K}^{\A}(\si(a))$, we also have 
$tp_{\langle K, \bar C_{K}\rangle}^{\A}(a)=tp_{\langle K, \bar C_{K}\rangle}^{\A}(\si(a))$. 
\par By Proposition \ref{Tsn}, $tp_{K}^{\A}(a)$ is isolated. So, by Lemma \ref{type}, it implies that $tp_{\langle K, C_{\A}\rangle}^{\A}(a)=tp_{\langle K, C_{\A}\rangle}^{\A}(\si(a)).$ Thus the map $\tilde \si$ sending $a$ to $\si(a)$ and which is the identity on $\langle K, C_{\A}\rangle$ is elementary. Going to a $\vert C_{\A}\vert^+$-saturated extension $\tilde \A$ of $\A$, we can extend it to an automorphism of that extension and 
since $L$ is strongly normal, 
we have $\langle L,C_{\tilde \A}\rangle=\langle K,\si(a), C_{\tilde \A}\rangle$. By Remark \ref{ldA}, $\langle L,C_{\tilde \A}\rangle\cap \U=\langle L,C_{\U}\rangle$ and $\langle \si(L),C_{\tilde \A}\rangle\cap \U=\langle \si(L),C_{\U}\rangle$. So, $\tilde \si\in Gal_{\U}(L/K)$. The uniqueness is clear. 
\par In the case when $\U=L$, we get that $gal(L/K)$ embeds in $Gal(L/K)$. 
\qed
\medskip

\lem\label{ann} Let $F\subset F_{0}\subset F_{1}$ be a tower of fields, $\bar u\in F_{1}$ 
and $q(\bar X,\bar Y)\in \IZ[\bar X,\bar Y]$, where $\bar Y:=(Y_{1},\cdots,Y_{n})$. 
\par Then the set of $n$-tuples $\bar c\in F_{0}^n$ such that $q(\bar X,\bar c)\in \I_{F_{0}}(\bar u)$ is a quantifier-free $F_{0}$-definable subset of $F_{0}^{n}$. 
\par Moreover, if $F(\bar u)$ and $F_{0}$ are linearly disjoint over $F$, then the set of $n$-tuples $\bar c\in F_{0}^n$ such that $q(\bar X,\bar c)\in \I_{F_{0}}(\bar u)$ is a quantifier-free $F$-definable subset of $F_{0}^{n}$.
\elem
\pr Let $\I:=\I_{F_{0}}(\bar u)$.
Since the polynomial ring $F_{0}[\bar X]$ is Noetherian, 
this ideal is finitely generated by say $f_{1},\cdots, f_{e}\in F_{0}[\bar X]$.  
Let us rewrite these generators $f_{m}(\bar X)$ as $g_{m}(\bar X,\bar c_{m})$, with $g_{m}(\bar X,\bar Z _{m})\in \IZ[\bar X,\bar Z_{m}]$ and $\bar c_{m}\in F_{0}$, $1\leq m\leq e$. 
\cl Let $p(\bar X)\in F_{0}[\bar X]$. Then we have $$(\bar F_{1}\models \forall \bar x(\bigwedge_{m=1}^e f_{m}(\bar x)=0\rightarrow p(\bar x)=0))\,\;{\rm  iff}\;\; p(\bar X)\in \I.$$
\ecl
\prc $(\Leftarrow)$ If $p\in \I$, then there exist $p_{m}(\bar X)\in F_{0}[\bar X]$, $1\leq m\leq e$, such that $p(\bar X)=\sum_{m=1}^e f_{m}(\bar X).p_{m}(\bar X)$ and so the implication holds. 
\par $(\Rightarrow)$ If $\bar F_{1}\models \forall x\;(\bigwedge_{m=1}^e f_{m}(\bar x)=0\rightarrow p(\bar x)=0)$, then since by hypothesis the $f_{m}'s$ belong to $\I$, we have $\bigwedge_{m=1}^e f_{m}(\bar u)=0$. So, we have $p(\bar u)=0$ and $p(\bar X)\in \I$. \qed
\medskip 
\par Now we rewrite $f_{m}(\bar X)$ as $g_{m}(\bar X,\bar c_{m})$, with $\bar c_{m}\in F_{0}$. Let $\bar Z:=(\bar Z_{m})_{1\leq m\leq e}$ and consider each $g_{m}(\bar X,\bar Z_{m})$ as an element of $\IZ[\bar X,\bar Z]$.
\par Let $\theta(\bar z, \bar y)$ be a quantifier-free $\L_{rings}$-formula equivalent in $ACF$ to the formula \\$\forall \bar x(\bigwedge_{m=1}^e g_{m}(\bar x,\bar z)=0\rightarrow q(\bar x,\bar y)=0),$ with $q(\bar X,\bar Y)\in \IZ[\bar X,\bar Y]$.
\par Let $\bar c\in F_{0}^{n}$. We have $F_{0}\models\theta((\bar c_{m})_{1\leq m\leq e},\bar c)$ if and only if it holds in $\bar F_{1}$. Since $\bar F_{1}\models ACF$, this is equivalent to $\bar F_{1}\models  \forall x\;(\bigwedge_{m=1}^e f_{m}(\bar x)=0\rightarrow q(\bar x,\bar c)=0)$, and by the Claim, equivalent to $q(\bar X,\bar c)\in \I$. Set $\tilde \theta(\bar y):=\theta((\bar c_{m})_{1\leq m\leq e},\bar y)$. Then $q(\bar X,\bar c)\in \I$ if and only if $\tilde \theta(\bar c)$ holds in $F_{0}$.
\par Now assume that $F(\bar u)$ and $F_{0}$ are linearly disjoint over $F$, so $\I_{F_{0}}(\bar u)$ is of the form $\I_{F}(\bar u).F_{0}[\bar X]$ (\cite[Chapter 3, section 2, Corollary 1]{Lang}). Therefore we may assume that $f_{1},\cdots, f_{e}\in F[\bar X]$. So we obtain for $\bar c\in F_{0}$, that $q(\bar X,\bar c)\in \I$ if and only if $\tilde \theta(\bar c)$ holds in $F_{0}$, where now the parameters occurring in the formula $\tilde\theta$ belong to $F$.
\qed
\thm\label{Gal} Let $L/K$ be a strongly extension. 
Then $Gal(L/K)$ is isomorphic to a group interpretable in $C_{\A}$ 
and the interpretation is given by $\L_{rings}(C_{K})$-quantifier-free formulas.
\ethm
\pr We proceed as in \cite{PS} but in order to check that indeed we do not need the hypothesis that $C_{\A}$ is algebraically closed, we will review the steps of the proof. We will use Proposition \ref{Tsn} and its proof. We have that $L=K\langle a \rangle$ for some element $a$ constrained over $K$. 
Let $\xi(y)$ be a (quantifier-free) $\L_{rings,D}(K)-$formula isolating $tp_{K}^{\A}(a)$.
\cl\label{def} {\rm \cite{PS}} There is a bijection between the $\L_{rings,D}(K)$-definable subset of $\langle L,C_{\A}\rangle$: $G_{\A}:=\{b\in \langle L, C_{\A}\rangle:\;\langle L,C_{\A}\rangle\models \xi(b)\}=\{b\in \A:\;\A\models \xi(b)\}$ and $Gal(L/K)$.
\ecl
\prc Let $b\in G_{\A}$ and let $\si\in Hom_{K}(L,\A)$ sending $a$ to $b$. By Remark \ref{ldA}, $\langle L,C_{\A}\rangle\cap L=C_{L}$ and so extending $\si$ on $C_{\A}$ by the identity and using that $L/K$ is strongly normal, we may extend $\si$ to an element of $Gal(L/K)$.
\par Conversely, let $\si\in Gal(L/K)$. Then $\si$ is uniquely determined by $\si(a)$ and the fact that it is the identity on $C_{\A}$. \qed
\medskip
\cl \label{G} {\rm \cite{PS}} One can define a group law $*$ on $G_{\A}$ by an $\L_{D}(K\cup\{a\})$-formula such that $(G_{\A},*)$ is isomorphic to the group $(Gal(L/K),\circ)$, where $\circ$ denotes the composition of maps. 
\ecl
\prc Let $\si,\;\tau\in Gal(L/K)$ and let $\mu:=\tau\circ\si$ with $\si(a)=b_{1},\;\tau(a)=b_{2},\;\mu(a)=b_{3}$. So $\tau(b_{1})=b_{3}$.
We have that $b_{1}=\tilde f_{\ell}(a, \bar k, \bar \alpha)$, with $\tilde f_{\ell}(Y,\bar T,\bar Z)=\frac{s_{1\ell}(Y,\bar T,\bar Z)}{s_{2\ell}(Y,\bar T,\bar Z)}$ the $\L_{fields,D}$-term obtained in Remark \ref{uni}, $\bar k\in K$, $\bar \alpha:=(\bar \alpha_{1}, \bar \alpha_{2})\in C_{\A}$, and $\ell$ vary over a finite set $\Lambda$ only depending on $\xi$. So, $\tau(b_{1})=b_{3}=\tilde f_{\ell}(\tau(a),\bar k, \bar \alpha)=\tilde f_{\ell}(b_{2},\bar k,\bar \alpha).$ This allows us to define the graph of a group law $*$ as a ternary relation $R$ between $b_{1},\; b_{2}$ and $b_{3}$ in $L$, as follows:
$R(b_{1}, b_{2}, b_{3}) \leftrightarrow \xi(b_{1})\;\&\;\xi(b_{2})\;\&\;\xi(b_{3})\;\&$
\begin{equation}
\bigvee_{\ell\in \Lambda}\exists \bar \alpha_{1}\bar \alpha_{2}\;(\bigwedge_{j=1}^2 D(\bar \alpha_{j})=0\;\&\;b_{1}=\tilde f_{\ell}(a, \bar k,\bar \alpha)\;\&\;b_{3}=\tilde f_{\ell}(b_{2},\bar k,\bar \alpha)).
\end{equation}
Then we express that this ternary relation is the graph of a group law. Note that this is existentially $K\cup\{a\}$-definable in $\langle L,C_{\A}\rangle$.
\medskip
\par Now we want to show that we can interpret this group in $C_{\A}$, using only parameters from $C_{K}$. At this point, we need to proceed in a more explicit way than in \cite{PS} since $C_{K}$ is not necessarily algebraically closed.
\par Let $b\in \langle L,C_{\A}\rangle$ such that $\langle L,C_{\A}\rangle\models \xi(b)$. Since $DCF_{0}$ admits q.e., we may assume that $\xi(y)$ is of the form $\bigwedge_{m} g_{m}(y)=0\;\&\;h(y)\neq 0$, $g_{m}(Y),\; h(Y)\in K\{Y\}$.
Substituting in that formula $\sigma(a):=\tilde f_{\ell}(a, \bar k, \bar \alpha)$ for $y$, where $b=\sigma(a)$, 
we obtain that some polynomials are in the annihilator of $(a^{\n},\bar k)$ in $C_{\A}[Y^*,\bar T]$ and that some (and in particular $s_{2\ell}^*(Y^*, \bar T, \bar \alpha_{2})$) do not belong to that annihilator.
 \par We rewrite the rational differential functions corresponding to $g_{m}(\tilde f_{\ell}(a,\bar k,\bar \alpha))$ as $\frac{r_{m1\ell}(a^{\n},\bar k,\bar \alpha_{1})}{r_{m2\ell}(a^{\n},\bar k,\bar \alpha_{2})}$ (respectively $h(\tilde f_{\ell}(a,\bar k,\bar \alpha))$ as $\frac{w_{1\ell}(a^{\n},\bar k,\bar \alpha_{1})}{w_{2\ell}(a^{\n},\bar k,\bar \alpha_{2})}$), where $r_{m1\ell}, r_{m2\ell}, w_{1\ell}, w_{2\ell}$ are $\L_{rings}$-terms.
 Note that $r_{m2j}(a^{\n},\bar k,\bar \alpha_{2}), \;w_{2\ell}(a^{\n},\bar k,\bar \alpha_{2})$ 
 are powers of $s_{2\ell}^*(a^{\n},\bar k,\bar \alpha_{2})$. 
 We will consider these terms as polynomials in $\IZ[Y^*,\bar T,\bar Z]$, with $Y^*$  of the same length as $a^{\n}$,
  $\bar T$ of the same length as the parameters $\bar k$ and $\bar Z$ of the same length as the constants $\bar \alpha:=(\bar \alpha_{1}, \bar \alpha_{2})$ (and assume the length of that last tuple is $d$). 
\par Let $\I:=\I_{C_{A}}(a^{\n},\bar k)$. By Lemma \ref{ann} (with $F_{0}=C_{\A}$ and $F_{1}=\langle L,C_{\A}\rangle$), the set of tuples $\bar c=(\bar c_{1},\bar c_{2})\in C_{\A}^{d}$ such that $(\bigwedge_{m}\;r_{m1\ell}(Y^*,\bar T,\bar c_{1})\in \I$ and \\$\;w_{1i\ell}(Y^*,\bar T,\bar c_{1}).s_{2\ell}^*(Y^*, \bar T,\bar c_{2})\notin \I)$
is definable by a quantifier-free $\La_{rings}$-formula $\theta_{\ell}(\bar t)$ with parameters in $C_{\A}$ 
 such that $C_{\A}\models \bigvee_{\ell\in \Lambda}\;\theta_{\ell}(\bar \alpha)$ if and only if $\langle L,C_{\A}\rangle\models \xi(\sigma(a))$. Set $\theta(\bar t):=\bigvee_{\ell\in \Lambda}\;\theta_{\ell}(\bar t)$.
 \par Since $L$ and $C_{\A}$ are linearly disjoint over $C_{K}$ (Remark \ref{ldA}), we have by Lemma \ref{ann}
 that $\I_{C_{\A}}(a^{\n},\bar k)=\I_{C_{K}}(a^{\n},\bar k).C_{\A}[\bar X].$ 
 Therefore we may assume that the parameters occurring in $\theta$ belong to $C_{K}$.
\medskip
\par Let $S:=\{\bar c:=(\bar c_{1},\bar c_{2})\in C_{\A}^{d}: C_{\A}\models \theta(\bar c)\}$.
\par Suppose that $\bar c\in S$, then for some $\ell\in \Lambda$, 
\begin{equation}
\bigwedge_{m}\;(r_{m1\ell}\in \I\;\&\;w_{1\ell}.s_{2\ell}^*\notin \I),
\end{equation}
 and the map sending $a$ to $\;b:=\tilde f_{\ell}(a,\bar k,\bar c)\in L$ is a partial elementary map. Indeed, $\bigwedge_{m}\;g_{m}(\tilde f_{j}(a,\bar k,\bar c))=0$ and $h(\tilde f_{j}(a,\bar k,\bar c))\neq 0$.
\medskip
\par We define a map $\nu$ from $C_{\A}^{d}$ to $\langle L,C_{\A}\rangle$, sending $\bar c\in S$ to one of the tuple $(\tilde f_{\ell}(a,\bar k,\bar c)\in G_{\A},$ where $\ell\in \Lambda$ is such that $(2)$ holds. (We may assume that $\Lambda$ is included in $\IN$ and we choose the smallest such index $\ell$.) Note that $\nu$ is definable.
\cl \label{nu} The following relation on $S$: $$\bar c\sim \bar c'\;\;{\rm iff}\;\;\nu(\bar c)=\nu(\bar c').$$
 is a $C_{K}$-definable equivalence relation in $C_{\A}$.
\ecl
\par $\bar c\sim\bar c'$ iff $\bigvee_{\ell\neq \ell'\in  \Lambda} (\theta_{\ell}(\bar c)\;\&\;\theta_{\ell'}(\bar c')\;\&\;\bigwedge_{j<\ell}\neg\theta_{j}(\bar c)\;\&\;\bigwedge_{j<\ell'}\neg\theta_{j}(\bar c')\;\&\;$
$$\frac{s_{1\ell}(a, \bar k, \bar c_{1})}{s_{2\ell}(a,\bar k,\bar c_{2})}=\frac{s_{1\ell}(a,\bar k,\bar c_{1}')}{s_{2\ell'}(a,\bar k,\bar c_{2}')})$$
$\;\;\;\;\;\;\;\;\;\;\;\;\;$ iff 
$\bigvee_{\ell\neq \ell'\in  \Lambda}( \theta_{\ell}(\bar c)\;\&\;\theta_{\ell'}(\bar c')\;\&\;\bigwedge_{j<\ell}\neg\theta_{j}(\bar c)\;\&\;\bigwedge_{j<\ell'}\neg\theta_{j}(\bar c')\;\&$
$$\;s^*_{1\ell}.s^*_{2\ell'}(\bar Y^*,\bar T,\bar c_{1},\bar c_{2}')-s^*_{1\ell'}.s^*_{2\ell}(\bar Y^*,\bar T, \bar c_{1}',\bar c_{2})\in \I).$$
\par Again we apply Lemma \ref{ann} and the fact that $L$ and $C_{\A}$ are linearly disjoint over $C_{K}$, and we get that this is a quantifier-free $\La_{rings}(C_{K})$-definable relation in $C_{\A}$. 
\qed
\medskip
\par 
\cl \label{R} $(S/\sim,R/\sim)$ is a group $\La_{rings}(C_{K})$-quantifier-free interpretable in $C_{\A}$.
\ecl
\par First we note that the relation $\sim$ is a congruence on $(S,R)$. Indeed, we have that if $\bar c\sim \bar d, \bar c'\sim \bar d', \bar c''\sim \bar d''$ and $R(\nu(\bar c), \nu(\bar c'), \nu(\bar c''))$, then $R(\nu(\bar d), \nu(\bar d'), \nu(\bar d''))$.
\par It remains to show that the relation $R$ is interpretable in $C_{\A}$. In equation $(1)$, assume that $b_{1}=\nu(\bar c_{11},\bar c_{12}),\;b_{2}=\nu(\bar c_{21},\bar c_{22}),\;b_{3}=\nu(\bar c_{31},\bar c_{32})$, namely 
\begin{equation}  
 b_{1}=\tilde f_{\ell}(a, \bar k,\bar c_{1}),\;b_{2}=\tilde f_{\ell}(a, \bar k,\bar c_{2}),\;b_{3}=\tilde f_{\ell}(a, \bar k,\bar c_{3}). 
\end{equation} 
Then express
$\;b_{3}=\tilde f_{\ell}(b_{2},\bar k,\bar c_{1})$
by substituting for $b_{3}$ and $b_{2}$ the expressions corresponding to equation $(3)$, and use Lemma \ref{ann} to get rid of the parameters $a^{\n},\;\bar k$.
\qed
\medskip
\par
\cl The group $(Gal(L/K),\circ)$ is isomorphic to $(S/\sim,R/\sim)$. 
\ecl
\par We send $\si\in Gal(L/K)$ to $\bar c/\sim$, where $\bar c\in C_{\A}^{d}$ is such that $\nu(\bar c)=\sigma(a)$.
\qed
\medskip
\par Using the same notations as in the proof above, we obtain the following corollary.
\cor\label{GaloisU} Let $L\subset \U\subset \A$ with $L/K$ strongly normal. \\Let $G_{\U}:=\{b\in \langle L, C_{\U}\rangle:\;\langle L,C_{\U}\rangle\models \xi(b)\}=\{b\in \U:\;\U\models \xi(b)\}$. 
Then there is a bijection between $Gal_{\U}(L/K)$ and $G_{\U}$.
\ecor
\pr Let $b\in G_{\U}$ and let $\si\in Hom_{K}(L,\U)$ sending $a$ to $b$. By Lemma \ref{emb}, $\si$ extends to an element of $Gal_{\U}(L/K)$.
\par Conversely, let $\si\in Gal_{\U}(L/K)$. Then $\si$ is uniquely determined by $\si(a)$ and the fact that it is the identity on $\langle K, C_{\U}\rangle$. Since $\si$ is the identity on $K$, $\A\models \xi(\si(a))$ and this quantifier-free formula also holds in $\U$ and $\langle L, C_{\U}\rangle$.
\qed
\medskip
\par So we obtain the following result which follows from former work of Kolchin {\rm \cite{K73}}.
\thm \label{Kolchin}  Let $L/K$ be a strongly extension. Then, the group $Gal(L/K)$ is isomorphic to the $C_{\A}$-points of an algebraic group $H$ defined over $C_K$.
Moreover for any intermediate subfield $L\subset \U\subset \A$, $Gal_{\U}(L/K)$ is isomorphic to $H(C_{\U})$. In particular, $gal(L/K)$ is isomorphic to $H(C_{K})$.
\ethm 
\pr The field of constants $C_{\A}$ is algebraically closed, so we apply the fact that $ACF_{0}$ admits e.i. So we have a $C_K$-definable function $f$ in $ACF_{0}$ which sends any equivalence class of $\sim$ on $C_{\A}^{n}$ to a unique element in $C_{\A}^m$.  By \cite[Proposition 3.2.14]{Ma}, there are finitely many $C_{K}$-algebraic subsets $X_{i}$ of  $C_{\A}^{n}$ such that $f\restriction X_{i}$ is a rational function with coefficients in $C_{K}$ (see section \ref{basicalg}). We already showed that the group law was definable over $C_K$. 
\par Therefore, if we take a tuple $\bar \alpha\in C_{\U}^{n}$ in an equivalence class of $\sim$, then $f(\bar \alpha)=\bar \beta\in C_{\U}^{m}$; furthermore the image of $\bar \alpha$ by the map $\nu$ (introduced above Claim \ref{nu}), is an element $b\in \langle L,C_{\U}\rangle\cap \xi(\U)$.
Then we apply Corollary \ref{GaloisU}.
\qed
\medskip
\nota \label{mu} We will denote by 
$\mu$ 
the isomorphism either from $Gal_{\U}(L/K)$ to $H(C_{\U})$ or its restriction from $gal(L/K)$ to $H(C_{K})$.
\enota
\cor \label{GalE} Let $L/K$ be a strongly extension and let $E$ an intermediate subfield.
Then $Gal_{\U}(L/E)$ is isomorphic to an $\L_{rings}(C_{K})$-quantifier-free definable subgroup of $H(C_{\U})$.
\par If $E$ is a strongly normal extension of $K$, then $Gal_{\U}(L/E)$ is a normal subgroup of $Gal_{\U}(L/K)$.
\ecor
\pr Note that $L$ is a strongly normal extension of $E$ and that $E$ is finitely generated over $K$, say $E:=K\langle \bar e\rangle$ with $\bar e\subset L$ ($L$ is also finitely generated as a field extension of $K$ and so one applies \cite[Chapter 3, section 2, Proposition 6]{Lang}). Again we may appeal to the result of Pogudin (recalled above Remark \ref{uni}) since $E\neq C_E$ is finitely generated and differentially algebraic and so w.l.o.g. $\bar e$ is reduced to one element $e$. 
\par We proceed as before noting that to express that $\si\in Gal_{\U}(L/E)$, it suffices to say that $\si\in Gal_{\U}(L/K)$ and $\si(e)=e$. So we express $e$ in terms of $a$ and $K$.
Namely,  $e=t(a, \bar k)$, with $t(Y,\bar T)\in \IQ\langle Y,\bar T\rangle$, $\bar k\in K$.
\par Assume that $\mu(\si):=\bar \alpha\in C_{\U}$, then, keeping the notation of Claim \ref{G}, we get: 
\begin{equation}
\si(a)=\tilde f(a, \bar k,\bar \alpha), \end{equation} 
\par Now we express that $\si(e)=e$ by:
\begin{equation}  e=t(\si(a), \bar k, \bar \alpha), \end{equation}
replacing $\si(a)$ by the expression $(4)$. 

So we obtain a $K\cup\{a\}$-definable subset $G_{E}$ of $G_{\U}$ such that $(G_{E},*_{T})$ is isomorphic to $(Gal_{\U}(L/E),\circ)$.
\par Then we replace $a$ by $a^{\n}$, similarly for $\bar k$ and the differential polynomials by their associated ordinary polynomials. We apply Lemma \ref{ann} as in Theorem \ref{Gal}. So $(5)$ is equivalent to express whether certain polynomials belong to $\I_{C_{\U}}(a^{\n}, \bar k^{\n})(=\I_{C_{K}}(a^{\n}, \bar k^{\n}).C_{\U}[\bar X])$. In this way, we express in the $\La_{rings}(C_{K})$-language which equivalence classes for the relation $\sim$ on $H(C_{\U})$ correspond to an element of $\mu(Gal_{\U}(L/E))$.
\par Let us prove the second assertion. Let $\tau\in Gal_{\U}(L/K)$, $\si\in Gal_{\U}(L/E)$ and $b\in E$. We want to show that $\tau(b)$ belongs to $\langle E, C_{\U}\rangle$. Since $E$ is a strongly normal extension of $K$ and $\tau\in Hom_{K}(E,\U)$, by Remark \ref{star}, $\langle E,C_{\U}\rangle=\langle \tau(E),C_{\U}\rangle$.
\qed
\medskip
\par Now we will examine under which hypothesis, a strongly normal extension of $K$ included in $\U$ is $Gal_{\U}(L/K)$-normal (see Definition \ref{normal}).
\prop \label{normalext} Assume that $\U$ is $\vert K\vert^+$-saturated. Suppose that $L/K$ is a strongly normal extension and that $dcl^{\U}(K)\cap L=K$. Then $L$ is a $Gal_{\U}(L/K)$-normal extension of $K$.
\eprop
\pr Let $L:=K\langle a\rangle$. By Lemma \ref{emb}, it suffices to prove that $L$ is $Hom_{K}(L,\U)$-normal.
\par Let $u\in L\setminus K$ and consider $tp_{K}^{\U}(u)$. Express $u$ as a rational function $\frac{p_{1}(a)}{p_{2}(a)}$, with $p_{1}(X),\;p_{2}(X)\in K\{X\}$ and $p_{2}(a)\neq 0$. By Proposition \ref{Tsn}, $tp_{K}^{\A}(a)$ is isolated by a quantifier-free formula $\chi(x)$ with parameters in $K$.
\par Suppose there is $d\neq u\in \U$ such that $tp_{K}^{\U}(u)=tp_{K}^{\U}(d)$. In particular for some $b\in \U$, we have $d=\frac{p_{1}(b)}{p_{2}(b)}\;\&\;\chi(b)$.
In particular, $tp_{K}^{\A}(b)=tp_{K}^{\A}(a)$. So, there is an element $\si\in Hom_{K}(L,\U)$ such that $\si(a)=b$.
\par So, either there is $u\neq d\in \U\setminus K$ such that $tp_{K}^{\U}(u)=tp_{K}^{\U}(d)$ and so we get that for some $\tilde\si\in Gal_{\U}(L/K)$, $\tilde\si(u)=d$.
\par Or $u$ is the only element of $\U$ realizing $tp_{K}^{\U}(u)$, then $u\in dcl^{\U}(K)$ and so by hypothesis $u\in K$. \qed
\medskip
\par Note that if $T_{c}=RCF$ and $dcl^{\U}(K)=K$, then $K\models T_{c}$.
 \medskip 
 \par From the above results we can deduce a partial Galois correspondence, between intermediate extensions and groups $\L_{rings}(C_{K})$-definable in $C_{\U}$.
\cor \label{interm} Let $K\subset L\subset \U$, suppose $L/K$ is strongly normal. Then,
\begin{enumerate}
\item the map sending $E$ to $Gal_{\U}(L/E)$ is an injective map from
  the set of intermediate extensions $E$ with $dcl^{\U}(E)\cap L=E$,
  to the set of groups  
$\L_{rings}(C_{K})$-definable in $C_{\U}$,
 \item assume that $K\subset E\subset L$ is a strongly normal extension of $K$, then $Gal_{\U}(L/K)/Gal_{\U}(L/E)$ is isomorphic to a subgroup of $Gal_{\U}(E/K)$ and the group $Gal_{\U}(L/K)/Gal_{\U}(L/E)$ is interpretable in $C_{\U}$.
 \end{enumerate} 
\ecor
\pr 
\par $(1)$ It follows from Proposition \ref{normalext} and Corollary \ref{GalE}.
\par $(2)$ Let $f$ be the map sending $\si\in Gal_{\U}(L/K)$ to $\si\restriction\langle E,C_{\U}\rangle$; it is well-defined and it goes from $Gal_{\U}(L/K)$ to $Gal_{\U}(E/K)$ by Remark \ref{star}. 
The kernel of $f$ is equal to $Gal_{\U}(L/E)$ (see also Corollary
\ref{Gal}).   
As in the proof of Theorem \ref{Kolchin}, we apply the elimination of imaginaries in $ACF_{0}$ and we obtain that the quotient $Gal(L/K)/Gal(L/E)$ is a group $H$ which is $C_{K}$-definable in $ACF_{0}$ and so an algebraic group defined over $C_{K}$. Since the definable function which chooses a point in each of the equivalence classes (containing a $C_{\U}$-point) is also $C_{K}$-definable, we obtain that $Gal_{\U}(L/K)/Gal_{\U}(L/E)$ is isomorphic to a subgroup of $H(C_{\U})$.
\par Since we have a Galois correspondance when dealing with the full galois group $Gal$, we have that $Gal(L/K)/Gal(L/E)$ is isomorphic to $Gal(E/K)$. Therefore, $Gal_{\U}(E/K)$ is isomorphic to the $C_{\U}$-points of an algebraic group isomorphic to $Gal(L/K)/Gal(L/E)$.\qed
\medskip
\par We will show that $Gal_{\U}(L/K)/Gal_{\U}(L/E)$ is, in general, isomorphic to a proper subgroup of $Gal_{\U}(E/K)$ (see Example \ref{proper}).

\section{Examples of strongly normal extensions}
\par We will provide in this section some examples of strongly normal extensions within some classes of differential fields and of their Galois groups. In particular, we review the two classical examples of Picard-Vessiot and Weierstrass extensions in classes of formally real fields or formally p-adic fields.  In this setting (and more generally), deep results on existence of strongly normal extensions have been obtained for instance in \cite{CHP}, \cite{GGO} and \cite{KP}. 
\par First we will recall the framework developed in \cite{GP}.
Then, we will review classical examples of strongly normal extensions as Picard-Vessiot extensions and Weierstrass extensions in these classes of topological differential $\L$-fields. 
\subsection{Model completion} Let $\La$ be a relational expansion of $\L_{rings}$ by $\{R_i; i\in
I\}\cup\{c_j;j\in J\}$ where the $c_j$'s are constants and the $R_i$'s 
are $n_{i}$-ary predicates, $n_{i}>0$. 
\dfn {\rm \cite{GP}}
Let $K$ be an $\La\cup \{^{-1}\}$-structure such that its restriction to $\La_{\text{rings}}\cup \{^{-1}\}$ is a field of characteristic $0$. Let $\tau$ be a  Hausdorff topology on $K$. Recall that $\langle K,\tau\rangle$ is a topological $\La$-field if the ring operations are continuous, the inverse function is continuous on $K\setminus\{0\}$ and every relation $R_{i}$ (respectively its complement $\lnot R_{i}$), with $i\in I$, is interpreted
in $K$ as the union of an open set $O_{R_{i}}$
(respectively $O_{\lnot R_i}$) and an algebraic subset $\{\bar x\in
K^{n_{i}}:\bigwedge_k r_{i,k}(\bar x)=0\}$ of $K^{n_i}$ (respectively
$\{\bar x\in K^{n_{i}}:\bigwedge_{l} s_{i,l}(\bar x)=0\}$ of $K^{n_{i}}$), where $r_{i,k},\;s_{i,l}\in K[X_1,\cdots,X_{n_i}]$.
\edfn
\par The topology $\tau$ is said to be first-order definable \cite{P87} if
there is a formula $\phi(x,\bar y)$ such that the set of subsets of the form $\phi(K,\bar a):=\{x\in K:\;K\models \phi(x,\bar a)\}$ with $\bar a\in K$ 
can be chosen as a basis $\V$ of neighbourhoods of $0$ in $K$. 
\par Examples of topological $\L$-fields are given in \cite[Section 2]{GP}. For instance, ordered fields, ordered valued fields, valued fields, $p$-valued fields, fields endowed with several distinct valuations or several distinct orders.
\par We now consider the class of models of a universal theory $T$ which has a model-completion $T_{c}$, and we assume that the models of  $T$ are in addition topological $\La$-fields, with a first-order definable topology.
Further, we work under the extra assumption that the class of models of $T_{c}$ satisfies {\it Hypothesis
$(I)$} \cite[Definition 2.21]{GP}. This hypothesis is the analog in our topological setting of the notion of {\it large fields} introduced by Pop \cite{Pop}.
\par We then consider the expansions of the models of $T$ to $\La_{D}$-structures and we denote by $T_D$ the $\La_D$-theory consisting of $T$ together the axioms expressing that $D$ is a derivation.
\par Hypothesis $(I)$ was used in \cite[Proposition 3.9]{GP} in order to show that any model of $T_D$ embeds in a model of $T_{c}$ which satisfies the scheme $(DL)$, namely that  if for 
each differential polynomial
 $f(X)\in K\{X\}\setminus\{0\}$, with $f(X)=f^*(X,X^{(1)},\ldots ,X^{(n)})$,
for every $W\in \V$, 

$\displaystyle (\exists \alpha_0,\ldots ,\alpha_n\;\;
f^*(\alpha_0,\ldots ,\alpha_n)=0 \wedge s_f^*(\alpha_0,\ldots ,\alpha_n)\ne 0 )
\Rightarrow$\\
$\displaystyle \Big(\exists z\big(f(z)=0\wedge s_f(z)\ne 0\wedge$
$\displaystyle \bigwedge_{i=0}^n (
z^{(i)}-\alpha_i\in W)\big)\Big)$.
\par Note that in \cite{GP}, the scheme $(DL)$ is not quite given as above, but in an equivalent form \cite[Proposition 3.14]{GP}.
\medskip
\par Note that since we assumed that the topology is first-order definable, the scheme
of axioms $(DL)$ can be expressed in a first-order way. Let $T_{c,D}^*$ be the $\La_{D}$-theory consisting of $T_{c}\cup T_{D}$ together with the scheme $(DL)$. Then $T_{c,D}^*$ is the model completion of $T_{D}$. A consequence of that axiomatisation is that the subfield of constants of a model $\U$ of $T_{c,D}^*$ is dense in $\U$ \cite[Corollary 3.13]{GP}.
\medskip
\nota Let $f\in K\{X\}$ and let $\langle f\rangle$ denote the differential ideal generated by $f$ in $K\{X\}$. Then, set $I(f):=\{g(X)\in K\{X\}: s_{f}^k.g\in \langle f\rangle \}.$ 
\par Let $\theta(\bar x)$ be a quantifier-free $\L_{D}$-formula. Then denote by $\theta^*(\bar y)$ the quantifier-free $\L$-formula obtained by replacing in $\theta(\bar x)$ every term of the form $x_{i}^{(j)}$ by a new variable $y_{ij}$ (\cite[Definition 3.16]{GP}).
\enota
\par In this section, $\U$ is a saturated model of $T_{c,D}^*$ containing $K$ and $\A$ as before a saturated model of $DCF_{0}$ containing $K$. 
(Recall that in a complete theory $T$, if the cardinal $\kappa$ is uncountable and strongly inaccessible, then there is a saturated model of $T$ of cardinality $\kappa$ \cite[Corollary 4.3.14]{Ma}. )
\par To stress that we consider strongly normal extensions within the class of models of $T_{D}$, we will use the term: {\it $T$-strongly normal}.
\dfn Let $K\models T_{D}$ and $L$ be a strongly normal extension of $K$. Then $L$ is a $T$-strongly normal extension of $K$ if $L$ is in addition a model of $T_{D}$.
\edfn
\lem Let $K\models T_{D}$. Let $F/K$ be a strongly normal extension, included in $\A$ and assume $F\neq K$; let 
$a\in \A$ be such that $F=K\langle a\rangle$. Let $\phi(x)$ be a formula isolating $tp^{DCF_{0}}_{K}(a)$. 
Then if there is $\bar b\in \U$ such that $\phi^*(\bar
b)$ holds, then $K$ has a $T$-strongly normal extension in $\U$, isomorphic to $F$ over $K$.
\elem
\pr  By \cite[Lemma 1.4]{M}, $\I^D_{K}(a)$ is of the form $I(f)$ for some $f\in K\{X\}$. Then we may assume that $\phi(x)$ is of the form $f(x)=0\;\&\;s_{f}.q(x)\neq 0$, where $q(X)\in K\{X\}$.
\par Since $f^*(\bar b)=0\;\&\;s_f^*(\bar b)\neq 0$, one can apply the scheme
$(DL)$. Therefore, one can find $u\in \U$ such that $f(u)=0$, with $u^{\n}$ close to $\bar b$. So, we also have that $s_f.q(u)\neq 0$. So $\phi(u)$ holds and by Proposition \ref{Tsn}, $K\langle u\rangle$ is a strongly normal extension of $K$ included now in $\U$ and so a model of $T_{D}$.
\qed
\medskip
\par Let $G$ be a connected algebraic group over $C_{K}$ and let $A\in LG(K)$, where $LG$ denotes the Lie algebra of $G$. If $C_{K}$ is existentially closed in $K$, then there is a strongly normal extension $L$ of $K$ for the logarithmic differential equation $dlog_{G}(Y)=A$ \cite[Theorem 1.3]{KP}. 
Moreover if $C_{K}$ is existentially closed in $K$, $C_{K}$ is large and is bounded (namely has for each $n$ only finitely many extensions of degree $n$, the so-called Serre property), then $K$ has a strongly normal extension $L$ in which $C_{K}$ is existentially closed \cite[Theorem 1.5]{KP}.
\par Further assume that $C_{K}\models T_{c}$, then both $K$ and the strongly normal extension $L$ of $K$ constructed above are models of $T_{D}.$  
\medskip
\par In Proposition \ref{normalext}, we assume that $dcl^{\U} K\cap L=L$, where $\U\models T_{c,D}$ and $\U$ is sufficiently saturated. We will denote by $\U\restriction \L$ the
$\L$-reduct of $\U$. In the following lemma, we
relate the algebraic closure in models of $T_{c,D}^*$ and of $T_{c}$.
\lem \label{acl} Let $L\models T_{D}$ and let $\U$ be a model of $T_{c,D}^*$
extending $L$. Then the algebraic closure $acl^{\U}(L)$ is equal to
$acl^{\U\restriction \L}(L)$.
\elem
\pr Let $a\in acl^{\U}(L)$ and let $\phi(x)$ be an
$\L_{D}(L)$-formula such that $\phi(a)$ holds in $\U$ and which has
only finitely many realizations. Since $T_{c,D}^*$ admits quantifier
elimination, $\phi(x)$ is equivalent to a finite disjunction of
formulas of the form: $$\bigwedge_{i\in I}
p_{i}(x)=0\;\&\;\theta(x),$$ where $\theta^*(\U)$ is an open subset of
some cartesian product of $\U$.
If $I=\emptyset$, then we obtain a contradiction since $\theta(\U)$ is infinite (it is a direct consequence of the scheme (DL) that near every tuple, one can find a tuple of the form $d^{\n}$ \cite[Lemma 3.12]{GP}).
Therefore we may assume that $I\neq \emptyset$; consider $\I^D_{L}(a)$. This is a prime
ideal of the form $\I(f)$, for some $f\in L\{X\}$ \cite[Lemma
1.4]{M}. Note that $f^*(a^{\n})=0\;\&\;s_f^*(a^{\n})\neq 0$. Either the set $S$ of solutions in $\U$ of the formula $f^*(\bar y)=0\;\&\;\theta^*(\bar y)$ is finite and so $a\in acl_{\L}(L)$. Or $S$ is infinite and for any $a\in S$, there exists a neighbourhood $V_{a}$ of $a$, included in $\theta^*(\U)$. Applying the scheme $(DL)$, there exist element $b\in\U\cap V_{a}$ satisfying $f(b)=0$ (and so $\theta(b)$ holds). So we get a contradiction with the finiteness of the number of solutions of $\phi(x)$.  
\qed
\medskip

 \subsection{Picard-Vessiot extensions}
 \par Let $K$ be a differential field which is a model of $T_{D}$ and assume that $C_{K}\models T_{c}$.
\par Let $P(Y):=Y^{(n)}+Y^{(n-1)}.a_{n-1}+\cdots+Y.a_{0}$,  be a differential linear polynomial of order $n$ with coefficients in $K$.
 \dfn \cite[Definition 3.2]{M} Let $L$ be an extension of $K$. Then $L$ is a Picard-Vessiot extension of $K$ corresponding to the linear differential equation $P(y)=0$ if
 \begin{enumerate}
 \item $L$ is generated by the solutions of $P(y)=0$ in $L$,
 \item $C_{L}=C_{K}$,
 \item $P(y)=0$ has $n$ solutions in $L$ which are linearly independent over $C_{K}$.
 \end{enumerate}
 \edfn
 \par Recall that if $L$ be a Picard-Vessiot extension of $K$, then $L$ is a strongly normal \cite[section 9]{M}.
 \par If $C_{K}$ is algebraically closed, there exists a Picard-Vessiot extension of $K$ corresponding to the equation $P(y)=0$, and such extension is unique up to $K$-isomorphism \cite[Theorems 3.4, 3.5]{Magid}. 
 \par In our particular setting, we have the following result.
 The corresponding algebraic equation is $Y_{n}+Y_{n-1}.a_{n-1}+\cdots+Y_{0}.a_{0}=0;$ it has $n$ linearly $C_{K}$-independent solutions: $u_{i}:=(u_{i0},\cdots, u_{in}),\;i=1\cdots,n$, equivalently the determinant of the corresponding $n\times n$-matrix formed by the first $n$ coordinates of the $u_{1},\cdots,u_{n}$ is non zero. Let $K\subset \U$ with $\U\models T_{c,D}^*$. By the scheme of axioms $(DL)$, we can find $n$ solutions: $v_{1},\cdots, v_{n}\in \U$ of the differential equation $P(y)=0$ with $(v_{i}, v_{i}^{(1)},\cdots,v_{i}^{(n-1)})$ in any chosen neighbourhood of $(u_{i 0},\cdots, u_{i n-1})$. Moreover by choosing a small enough neighbourhood, we can guarantee that the Wronskian of $v_{1},\cdots, v_{n}$ is non zero, equivalently $v_{1},\cdots, v_{n}$ are linearly independent over $C_{K}$ \cite[Lemma 4.1]{M}. However the field of constants of $K(v_{1},\cdots,v_{n})$ might be bigger than $C_{K}$. 
\rem Assume that $C_{K}$ is existentially closed in the class of models of $T$ and let $S$ be the full universal algebra solution algebra for $P$ \cite[Definition 2.12]{Magid}. 
Then a corollary of the proof of  \cite[Theorems 3.4, 3.5]{Magid} is that if there is a maximal differential ideal $I$ of $S$ such that $Frac(S/I)$, the fraction field of $S/I$, is a model of $T_{D }$, then there exists a Picard-Vessiot extension of $K$ corresponding to the equation $P(y)=0$ \cite[Corollary 1.18]{Magid} (which will be $T$-strongly normal).
\erem
\medskip
\lem {\rm  \cite[Theorem 5.5]{Kap}} Let $L$ be a Picard-Vessiot extension of $K$, then $gal(L/K)$ is isomorphic to the $C_{K}$-points of a  linear  algebraic group defined over $C_{K}$. \qed
\elem
\medskip
\prop {\rm (\cite[Proposition 3.9]{Magid})} Let $L$ be a differential extension of $K$ and suppose that 
\begin{enumerate}
\item $C_{K}=C_{L}$,
\item there is a finite-dimensional $C_{K}$-vector space $V$ such that $L=K\langle V\rangle$,
\item there is a subgroup $H$ of $gal(L/K)$, leaving $V$ invariant and such that $L/K$ is $H$-normal. 
\end{enumerate}
Then $L/K$ is a Picard-Vessiot extension of $K$. \qed
\eprop
\cor {\rm \label{PV-normal}}
Let $K$ be a differential field and let $L$ be a finite Galois extension of $K$ such that $C_{K}$ is algebraically closed in $L$. Then $L$ is a Picard-Vessiot extension of $K$.\qed
\ecor
\ex\label{proper}
When $T$ is the theory of ordered fields, we provide here an example of extensions $K\subset E\subset L$ satisfying the hypotheses of 
	Proposition \ref{interm} and such that $Gal_{\U}(L/K)/Gal_{\U}(L/E)$ is isomorphic to a proper subgroup of $Gal_{\U}(E/K)$.
	
	Let $K:=\R$ endowed with the trivial derivation,
	$E:=\R \langle t\rangle$ where $D(t)=t$,
	$L:=\R \langle u\rangle$ where $u^2=t$ (and so $D(u)=\frac{1}{2}u$),
	$L/K$ and $E/K$ are Picard-Vessiot and so are strongly normal.
	Take $\sigma\in Gal_{\U}(E/K)$ such that $\sigma(t)=-t$. 
	Since $t$ is a square in $L$ and $-t$ is not, $\sigma$ may clearly not be in the image of the map $f$ from the proof of Proposition \ref{interm}.
	
	Here we can easily describe $Gal_{\U}(E/K)$ and $Gal_{\U}(L/K)$. They are both isomorphic to $\mathbb G_m(C_\U)$. 
	
	However  $Gal_{\U}(L/E)$ is isomorphic to the subgroup $\{-1,1\}$ of $\mathbb G_m(C_\U)$. 
	Note that since $C_\U$ is real closed, any element of the group $\mathbb G_m(C_\U)/\{-1,1\}$ is a square, 
	which is no longer true for the group $\mathbb G_m(C_\U)$ and so they cannot be isomorphic.

\eex

\subsection{Weierstrass extensions}\label{W}
Now we will consider the case of a strongly normal extension generated by a solution of a differential equation $W_{k}(Y)=0$ where $W_{k}(Y):= (Y^{(1)})^2-k^2.(4.Y^3+g_{2}.Y-g_{3})\;(\star)$, with $k,\;g_{2},\;g_{3}\in K$ and $27g_{3}^2-g_{2}^3\neq 0$. 
Note that $s_{W_{k}}=\frac{\partial}{\partial Y^{(1)}} W_{k}=2.Y^{(1)}$.
Let $K$ be a differential field satisfying $T_{D}$ and let $\U$ be a model of $T_{c,D}^*$ containing $K$. Then given any solution $(u,u_{1})\in \U^2$ of the Weierstrass equation: $(\dagger)$ $Y^2=4.X^3+g_{2}.X-g_{3}$ with $u_{1}\neq 0$, there exists a solution $a$ of the differential equation $W_{k}(Y)=0$ such that $(a,\frac{a^{(1)}}{k})$ is close to $(u,u_{1})$. 
\par Recall that  an elliptic curve (over $K$) is defined as the projective closure of a nonsingular curve defined by the equation $(\dagger)$. We will denote by $\E$ the corresponding elliptic curve (or $\E(F)$ if we look at its points on an intermediate subfield $F$: $K\subset F\subset \A$) and usually we will use the affine coordinates of the points on $\E$.
We denote by $\oplus$ the group operation on the elliptic curve $\E$ and by $\ominus$ the operation of adding the inverse of an element. It is well-known that $(\E,\oplus)$ is an algebraic group over $K$ (one can express the sum of two elements as a rational function of the coordinates of each of them \cite[Chapter 3, section 2]{Sil}).
\par Let $\si\in Hom_{K}(L,\A)$ and let $a$ be such that $W_{k}(a)=0$. Then $\si(a)$ also satisfies $W_{k}(y)=0$. 
\par Consider the element $(\si(a),\frac{\si(a)^{(1)}}{k}) \ominus (a,\frac{a^{(1)}}{k})$, then one can verify that it belongs to $\E(C_{\A})$ (see \cite[Lemma 2, chapter 3]{K53}, or for instance \cite[section 9]{M}).
So we get the following Lemma.
\lem Let $L:=K\langle a\rangle$, where $a\in \U$ satisfies $W_{k}(a)=0$. Assume that $C_{K}=C_{L}$. Then $L$ is a  strongly normal extension of $K$.\qed
\elem
\par Assume that $L$ is a strongly normal extension of $K$ generated by $a$ with $W_{k}(a)=0$ and set $\mu: gal(L/K)\rightarrow \E(C_{K}): \si\rightarrow (\si(a),\frac{\si(a)^{(1)}}{k}) \ominus (a,\frac{a^{(1)}}{k})$ (respectively $\mu: Gal_{\U}(L/K)\rightarrow \E(C_{\U}): \si\rightarrow (\si(a),\frac{\si(a)^{(1)}}{k}) \ominus (a,\frac{a^{(1)}}{k})$). 
\medskip
\lem \rm{\cite[section 9, Example]{M}} The group $(gal(L/K),\circ)$ \\(respectively $(Gal_{\U}(L/K),\circ)$) is isomorphic to a definable subgroup of  $(\E(C_{K}),\oplus)$ (respectively $(\E(C_{\U}),\oplus)$).
\elem 
\pr It is easy to see that $\mu$ is an injective group morphism \cite[Chapter 3, section 6]{K53}. In order to show that the image of $\mu$ is a definable subset of $\E(C_{K})$ (respectively $\E(C_{\U})$), one notes that $(c_{1},c_{2})\in \E(C_{K})$ belongs to the image of $\mu$ if the first coordinate of $(c_{1},c_{2})\oplus(a,\frac{a^{(1)}}{k})$ satisfies the formula $\chi(x)$ isolating $tp_{K}^{DCF_{0}}(a)$. \qed
\medskip
 \par In order to identify $\mu(gal(L/K))$, we now examine in more details which kind of (definable) subgroups may occur in $\E(C_{K})$, when $C_{K}$ is first a real closed field and then a $p$-adically closed field.
\par Let $C_{K}\models RCF$. Then the o-minimal dimension of $\E(C_{K})$ assuming $\E(C_{K})$ infinite, is equal to $1$ and so in case $gal(L/K)$ is infinite, the o-minimal dimension of $\mu(gal(L/K))$ is also equal to $1$. It has a connected component of finite index $\mu(gal(L/K))^{0}$ which is definable \cite[Proposition 2.12]{P88}. So in case $C_{K}=\IR$, we get that $\mu(gal(L/K))^{0}$ is isomorphic to $\IR/\IZ$, as a connected $1$-dimensional compact commutative Lie group over $\IR$ .
\par Recall that over $\IR$, $\E(\IR)$ is either isomorphic to $S^1\times \IZ/2\IZ$ or to $S^1$, depending on the sign of the discriminant (or on whether the polynomial $4.X^3+g_{2}.X-g_{3})$ has three or one real roots) (see \cite{Sil2} or \cite{Poonen}). So the only non-trivial proper subgroup of finite index in $\E(\IR)$ has index $2$; this property (restricted to definable subgroups) can be expressed in a first-order way. So it transfers to any other real closed field, and so $\mu(gal(L/K))^{0}\cong \E(C_{K})^{0}$. Moreover we have that $\mu(gal(L/K))^{0}/\mu(gal(L/K))^{00}\cong \E(C_{K})^{0}/\E(C_{K})^{00}$ where $\mu(gal(L/K))^{00}$ is the smallest type-definable subgroup of bounded index of $\mu(gal(L/K)).$ 
\par Suppose now that $C_{K}$ is a $p$-adically closed field, then one can find a description of definable (and type-definable) subgroups of $\E(C_{K})$ in \cite[Chapter 7]{Sil} and \cite{OP}. Even though in this case the theory does not admit elimination of imaginaries, since we obtained the isomorphism explicitly, we get that $gal(L/K)$ (respectively $Gal_{\U}(L/K)$) is isomorphic to the $C_{K}$-points (respectively to its $C_{\U}$-points) of a group $\L_{rings}(C_{K})$-definable. However, even in the case $C_{K}=\IQ_{p}$, the lattice of subgroups is more complicated. Let us give a quick review of some of the known subgroups. \par First let us recall some notations. Given an elliptic curve $\E$, $\tilde \E$ is its reduction modulo $p$ and $\tilde \E_{ns}$ the subset of the non-singular points of $\tilde \E$. One can also attach to $\E$ a formal group $\hat{\E}$ over $\IZ_{p}$.
Then we have $\E_{1}(\IQ_{p})\subset \E_{0}(\IQ_{p})\subset \E(\IQ_{p})$ and $\E_{1}(\IQ_{p})$ has a subgroup of finite index isomorphic to the additive group of $\IZ_{p}$, where $\E_{0}(\IQ_{p}):=\{P\in \E(\IQ_{p}): \tilde P\in \tilde \E_{ns}(\IF_{p})\}$ and $\E_{1}(\IQ_{p}):=\{P\in \E(\IQ_{p}):\tilde P=\tilde 0\}\cong \hat{\E}(p.\IZ_{p})$ \cite[Proposition 6.3]{Sil}. The index $[\E(\IQ_{p}):\E_{0}(\IQ_{p})]=4$ and $\E_{0}(\IQ_{p})/\E_{1}(\IQ_{p})\cong \tilde \E_{ns}(\IF_{p})$.
\par In case $C_{K}$ is sufficiently saturated, we have $\E(C_{K})$ has an open definable subgroup $H$ such that $H/H^{\circ\circ}$ is isomorphic to the profinite group $(p^r \IZ_{p},+,0)$, where $r>\frac{v(p)}{p-1}$ \cite[Chapter 4, 6.4, b)]{Sil}.
\section{Ordered differential fields}
\par We will use the same notations as in the previous section. We will consider the case where $T$ is the theory of ordered fields and so $T_{c}$ is the theory $RCF$ of real closed fields and $T_{c,D}^*$ the theory $CODF$ of closed ordered differential fields. 
Note that in this case, the language $\La_{D}$ is equal to the language of ordered differential fields, namely $\La_{D}=\{+,-,.,{}^{-1},0,1,<,D\}$. We will use the fact that $CODF$ admits elimination of imaginaries \cite{Point}. 
\par Let $K\subset L\subset \U$ with $K, L\models T_{D}$. Assume that $C_K\subset_{ac} C_{\U}$ and that $L$ is a strongly normal extension of $K$, then 
both $gal(L/K)$ (respectively $Gal_{\U}(L/K)$) are isomorphic to the $C_{K}$-points (respectively the $C_{\U}$-points) of an algebraic group  $\L_{rings}(C_{K})$-definable. We will further assume that $\U\models CODF$ and so $C_{\U}\models RCF$. (Note that $C_K\subset_{ac} C_{\U}$ implies that $C_{K}\models RCF$, as well.)
\par By former results of A. Pillay \cite{P88}, any definable group $G$ in a real closed field can be endowed with a definable topology in such a way it becomes a topological group. Moreover $G$ has a connected component $G^0$ of finite index (which is also definable) \cite[Proposition 2.12]{P88}. If $G$ lives in a sufficiently saturated model, then it has a smallest type-definable subgroup of bounded index denoted by $G^{00}$; this result holds in fact in any $NIP$ theory \cite{Sh}. One can endow the quotient $G/G^{00}$ with the logic topology and if $G$ is in addition definably connected and definably compact, then $G/G^{00}$ is a compact connected Lie group \cite[Theorem 1.1]{BOPP}. Moreover the $o$-minimal dimension of $G$ is equal to the dimension $G/G^{00}$ as a Lie group \cite[Theorem 8.1]{HPP}.

\nota 
We denote by $aut(L/K)$ the group of $\L_{D}$-automorphisms of $L$ fixing $K$ and by $Aut_{T}(L/K)$ the group $aut(\langle L, C_{\U}\rangle/\langle K, C_{\U}\rangle)$.
\enota
\par Let $L^{rc}$ be the real closure of $L$ in $\U$. Note that any $\si\in aut(L/K)$ (respectively $aut(\langle L,C_{\U}\rangle/\langle K,C_{\U}\rangle)$ has a unique extension to an element of $aut(L^{rc}/K)$ (respectively $aut(\langle L^{rc},C_{\U}\rangle/\langle K,C_{\U}\rangle)$. In the next proposition
we will use the following easy consequence of Lemma \ref{acl} that
$L^{rc}=dcl^{\U}(L)$ (one uses in addition that in an ordered structure,
the algebraic closure is equal to the definable closure).
\par Let $L$ be a strongly normal extension of $K$ and $\bar d\in L^{rc}$. 
Then $L\langle \bar d\rangle=L(\bar d)$ and $L\langle\bar d\rangle$ is again a strongly normal extension of $K\langle\bar d\rangle$.
First $C_{L(\bar d)}=C_{L}=C_{K}$ since $C_{K}=C_{L}$ is algebraically closed in $L$. 
If $\si\in Gal_{\U}(L\langle\bar d\rangle/K\langle\bar d\rangle)$, $\si(L)\subset \langle L,C_{\U}\rangle$ and $\si(\bar d)=\bar d$. So, $\si(L(\bar d))\in \langle L(\bar d),C_{\U}\rangle$.
\par Since $\bar d\in L^{rc}=dcl^{\U}(L)$, any $\tau\in Aut_{T}(L/K)$ extends uniquely to an element of $Aut_{T}(L\langle\bar d\rangle/K\langle\bar d\rangle)$. In the proposition below, we will obtain a relative Galois correspondance for definable subgroups of $\mu(Aut_{T}(L/K))$ (see Notation \ref{mu}).
\prop\label{inter} Let $L/K$ be a strongly normal extension with $L\subset \U$.
\par Let $G_{0}$ be an $\L(K)$-definable subgroup of $H(C_{\U})$. Then there is a tuple of elements $\bar d\in L^{rc}$ such that we have: $G_{0}\cap \mu(Aut_{T}(L/K))\cong Aut(L(\bar d)/K\langle \bar d\rangle)$.
\eprop
\pr  Let $L:=K\langle a\rangle$ for some $a\in L$.
Let $G_{0}$ be a definable subgroup of $H(C_{\U})$ by some $\L(K)$-formula $\psi(\bar x)$. 

\par Consider $W:=\{\tau(a):\;\mu_{T}(\tau)\in G_{0}\;\wedge\;\tau\in Gal_{\U}(L/K)\}$, it is a $\L_{D}(L)$-definable subset of $\langle L,C_{\U}\rangle$. Indeed, we can write $W$ as $\{(\tilde f_{\ell}(a, \bar k, \bar \alpha))$, with $\tilde f_{\ell}(Y,\bar Z,\bar T)$ $\L_{D}\cup\{^{-1}\}$-terms as in Claim \ref{G}, with $\bar k\in K$, $\bar \alpha\in C_{\U}$ and $\ell$ vary over a finite set $\Lambda$.
Since $C_{\U}$ is definable in $\U$,
we can express in $\U$ that $C_{\U}\models\psi(\bar \alpha)$ by a $\L_{D}(K)$-formula .
\par Since $CODF$ has elimination of imaginaries \cite[Theorem 2.6]{Point}, there is a canonical parameter $\bar d$ for $W$ in some elementary extension $\tilde \U$ of $\U$, namely a tuple fixed by any $\L_{D}$-automorphism which leaves $W$ invariant \cite[Definition 8.4.1]{TZ}.  Note that $\bar d$ belongs to $dcl^{\tilde \U}(L)=L^{rc}$ since any $\L_{D}$-automorphism of $\tilde \U$ which fixes $L$,
leaves $W$ invariant.
\par Set $F:=K\langle \bar d\rangle$. Let us show that $G_{0}\cap \mu_{T}(Aut_{T}(L/ K))\cong Aut_{T}(L(\bar d)/F)$.
\par Let $\tau\in Aut_{T}(L( \bar d)/F)$ and consider $\tau\restriction L(\bar d)$. 
It extends to an automorphism $\tilde \tau$ of $\U$, which leaves $F$ fixed and so it leaves $W$ invariant.  Since $\bar a\in W$, $\tau(\bar a)\in W$, namely there is $\tau_{0}\in Gal_{\U}(L/K)$ with $\mu(\tau_{0})\in G_{0}$ and $\tau(\bar a)=\tau_{0}(\bar a)$.
Since $\mu(\tau_{0})$ is determined by the action of $\tau_{0}$ on $\bar a$, we get that $\mu(\tau)\in G_{0}$. 
\par Since $L$ is a strongly normal extension of $F$ and $\tau(L)\subset \U$, $\langle L,C_{\U}\rangle=\langle \tilde \tau(L),C_{\U}\rangle=\langle \tau(L),C_{\U}\rangle$ (see Remark \ref{star}). 
In particular, $\tilde\tau(a)=\tau(a)\in \langle L,C_{\U}\rangle$. 

\par Conversely, if $\bar c\in G_{0}\cap \mu(Aut_{T}(L/K))$, then $\mu^{-1}(\bar c)\in Aut_{T}(L/K)$.

 Again $\mu^{-1}(\bar c)\restriction L$ extends to an $\La_{D}$-automorphism $\si$ of $\U$, which has the same action on $\bar a$ and so leaves $W$ invariant and thus $\bar d$ fixed.
 Therefore the restriction of $\si$ to $\langle L(\bar d), C_{\U}\rangle$ belongs to $Aut_{T}(L(\bar d)/ F)$.
 \qed
\medskip
\par Note that we could have expressed the conclusion of the above proposition as follows.
$G_{0}\cap \mu_{T}(Aut_{T}(L^{rc}/K))\cong Aut_{T}(L^{rc}/K\langle \bar d\rangle)$.

\section{ On the non finitely generated case}\label{nfg}
J. Kovacic extended the notion of strongly normal extensions to the infinitely generated case \cite[Chapter 2, section 1]{Kov1}. Recall that $L$ is a strongly normal extension of $K$ if $L$ is a union of finitely generated strongly normal extensions of $K$. Equivalently, $C_{L}=C_{K}$ and $L$ is generated by finitely generated strongly normal extensions of $K$ \cite[Chapter 2, Section 1, Definition]{Kov1}. 
\par In this section we will place ourselves in the setting of section 3. We will fix a differential field $K$, we work inside an extension $\U$, and embed $\U$ in a saturated model $\A$ of $DCF_0$.
\par In \cite[Chapter 2, Theorem 1]{Kov1}, J. Kovacic proved that if $L$ is a strongly normal extension of $K$ with $C_{K}$ algebraically closed, then on one hand, the set of finitely generated strongly normal extensions of $K$ in $L$ forms an injective system and on the other hand, $Gal(L/K)$ is the inverse limit of the differential Galois groups of finitely generated extensions. In particular, it is isomorphic to an inverse limit of algebraic groups.
\par In the following proposition, we will establish an analog for the differential Galois groups $gal(L/K)$ and $Gal_{\U}(L/K)$ (Notation \ref{auto-nota}) (without the assumption that $C_K$ is algebraically closed). 

\medskip
\prop Let $L/K$ be a strongly normal extension with $L\subset \U$.
Then $Gal_{\U}(L/K)$) is isomorphic to a subgroup of an inverse limit of groups of the form $H_i(C_{\U})$ where $H_i$ is an algebraic group defined over $C_{K}$, $i\in I$. Moreover, $gal(L/K)$ is a subgroup of $Gal_{\U}(L/K)$ and is isomorphic to a subgroup of an inverse limit of groups of the form $H_i(C_K)$, $i\in I$.
\eprop 
\pr Let $L=\bigcup_{F\in \F} F$, where $\F$ denotes the set of all finitely generated strongly normal extensions of $K$ in $L$. Then $Gal_{\U}(L/F)$ is a normal subgroup of $Gal_{\U}(L/K)$ (see Corollary \ref{GalE}, note that there we only use that $F$ was a (finitely generated) strongly normal extension of $K$). 
Moreover we have an embedding of $Gal_{\U}(L/K)/Gal_{\U}(L/F)$ into $Gal_{\U}(F/K)$, sending $\si$ to $\si\restriction F$ (Corollary \ref{interm}).

\par  For each $K\subset F_{i}\subset F_{j}\subset L$ with $F_{i},\,F_{j}\in \F$, define the maps $$f_{F_{i}}^{F_{j}}:Gal_{\U}(L/K)/Gal_{\U}(L/F_{j})\rightarrow Gal_{\U}(L/K)/Gal_{\U}(L/F_{i}),$$ sending $\si.Gal_{\U}(L/F_{j})$ to $\si.Gal_{\U}(L/F_{i})$ and $$r_{F_{i}}^{F_{j}}:Gal_{\U}(F_{j}/L)\rightarrow Gal_{\U}(F_{i}/K)$$ sending $\si$ to $\si\restriction F_{i}$ (with kernel $Gal_{\U}(F_{j}/F_{i})$).
The following diagram commutes:
\medskip

\begin{center}
$\begin{CD}
Gal_{\U}(L/K)/Gal_{\U}(L/F_{j}) @>>> Gal_{\U}(F_{j}/L)\\
@Vf_{F_{i}}^{F_{j}}VV  @Vr_{F_{i}}^{F_{j}}VV\\
Gal_{\U}(L/K)/Gal_{\U}(L/F_{i}) @>>> Gal_{\U}(F_{i}/K)
\end{CD}$
\end{center}
\medskip
\par The family of differential Galois groups $Gal_{\U}(F/K)$ forms a projective system, as well as the family of quotients $Gal_{\U}(L/K)/Gal_{\U}(L/F)$, $F\in \F$. \\Since $\bigcap_{F\in \F} Gal_{\U}(F/K)=\{1\}$, we get an embedding of $Gal_{\U}(L/K)$ into $\varprojlim_{F\in \F} Gal_{\U}(F/K)$.
\par Finally, for $\si\in Gal_{\U}(L/K)$, define $\mu(\si):=(\mu(\si\restriction F))$ (Notation \ref{mu}).\\
We get $\mu(Gal_{\U}(L/K))\subset \varprojlim_{F\in \F}\mu(Gal_{\U}(F/K))$.
\par The second part is clear. \qed
\medskip
\par Let us give an example of a strongly normal extension which is possibly infinitely generated. 
\ex Let $K$ be a formally real differential field with a real closed subfield of constants.  Let $K^*$ be the intersection of all real closed subfields of a given algebraic closure of $K$. Then by a result of M. Griffin \cite{C}, $K^*$ is the maximal Galois extension of $K$ to which all the orderings of $K$ extend \cite[Theorem 1.1]{C}. Then by Corollary \ref{PV-normal}, $K^*$ is a union of Picard-Vessiot extensions of $K$. 
\eex
\par Adapting the classical proofs \cite{Magid}, we will show that one can find a maximal strongly normal extension of $K$ inside $\U$ which has proper strongly normal extension (inside $\U$).
\rem Let $K\models T_{c}$ with a trivial derivation and suppose that $K$ is $\aleph_{0}$-saturated (we have to be able to find $n\geq 3$ elements $\alpha_{1},\cdots,\alpha_{n}$ of $K$ linearly independent over $\IQ$). Let $\U$ be a saturated extension of $K$ satisfying $T_{c,D}^*$. Inside $\U$, choose $n$ elements $a_{1},\cdots, a_{n}$ which are algebraically independent over $K$. Consider $L_{2}$ the differential subfield generated by $K$ and $a_{1},\cdots, a_{n}$. Consider the subfield $L_{1}$ of $L_2$ containing $K$ and generated by $b_{1},\cdots, b_{n}$ the coefficients of the polynomial $\prod_{i=1}^n (X-a_{i})$ (namely the $b_{i}$'s are the images (up to a sign) of the elementary symmetric functions on $a_{1},\cdots, a_{n}$ and so also algebraically independent over $K$). 
Now define a derivation on $L_{2}$ by first imposing that $D(b_{i})=\alpha_{i}.b_{i}$ and then $D(a_{i})$ is given by using the fact that $a_{i}$ is algebraic over $L_{1}$.
\par The field $L_{1}$ is a Picard-Vessiot extension of $K$
\cite[Proposition 1.20, Example 3.23]{Magid} corresponding to the linear differential
equation $Y^{(n)}+Y^{(n-1)}.\beta_{n-1}+\cdots+Y.\beta_{1}=0$,
where $\prod_{i=1}^n
(Y-\alpha_{i})=Y^n+Y^{n-1}.\beta_{n-1}+\cdots+Y.\beta_{1}$. The extension
$L_2/L_1$ is a Picard-Vessiot extension. Indeed, we use Corollary \ref{PV-normal}: $L_2$ is a finite
Galois extension of $L_{1}$ and $C_{L_{1}}=C_{K}=K$ is algebraically
closed in $\U$. 
However $L_{2}$ is not contained in a Picard-Vessiot extension of $K$. Suppose otherwise,
then $\langle L_{2},\bar C_{K}\rangle$ would also be contained in a Picard-Vessiot extension of $\langle K, \bar C_{K}\rangle$, which contradicts \cite[Example 3.33]{Magid}.
\par Now let $F$ be a maximal strongly normal extension of $K$ inside $\U$ containing $L_{1}$. Let us show that the subfield $F_{1}:=\langle F, L_{2}\rangle$ of $\U$ is a proper strongly normal extension of $F$. 
\par Indeed, $F_{1}/F$ is a Picard-Vessiot extension since it is a finite
Galois extension and $C_{F_1}=C_K$. Suppose that $F_1=F$, then
$\langle L_{2},\bar C_{K}\rangle$ is
contained in a strongly normal extension of $\langle K, \bar
C_{K}\rangle$. It implies by \cite[Proposition 11.4]{Kov3} that $\langle L_{2},\bar C_{K}\rangle$ is 
contained in a Picard-Vessiot extension of $\langle K,\bar
C_{K}\rangle$, a contradiction.
\erem
\bigskip
\noindent {\bf Acknowledgments:}
Both authors thank the MSRI for their hospitality during the spring of 2014 and the FRS-FNRS for their support. 
 They thank Omar Sanchez for his kind and useful remarks and an anonymous referee, which led to a significant improvement of the paper. They thank Anand Pillay for his insightful remarks, in particular for insisting on the fact that a previous hypothesis on the field of constants of the intermediate field $\U$ was not necessary, and on Lemma 4.4 (and the comment afterwards). We are also grateful to Michael Singer for his advice.

\end{document}